\documentclass[10pt,a4paper,reqno]{amsart}
\usepackage{amsthm}
\usepackage{amsmath}
\usepackage{amssymb}
\usepackage{graphicx,color}

\usepackage[font=small]{caption}
\usepackage{subcaption}

\usepackage{multirow}
\usepackage[shortlabels]{enumitem}

\usepackage{soul} 

\makeatletter
\theoremstyle{plain}
\newtheorem{theorem}{Theorem}
  \theoremstyle{definition}
  \theoremstyle{remark}
  \newtheorem{remark}[theorem]{Remark}
  \theoremstyle{plain}
  \newtheorem{Proposition}[theorem]{Proposition}
  \theoremstyle{plain}
  \newtheorem{lemma}[theorem]{Lemma}
  \theoremstyle{plain}
  
 \theoremstyle{definition}
  
\usepackage{amsfonts}
\usepackage{mathrsfs}

\newtheorem*{question*}{\it{QUESTION}}

\theoremstyle{plain}

\newtheorem*{oq*}{Open Question}

\newcommand{\N}{\mathbb{N}}
\newcommand{\R}{{\mathbb{R}}}
\newcommand{\C}{{\mathbb{C}}}

\newcommand{\dd}{{\rm d}}
\newcommand{\ii}{{\rm i}}
\newcommand{\e}{{\rm e}}
\newcommand{\ess}{{\rm ess}}

\newcommand{\pd}[2]{\frac{\partial #1}{\partial #2}}

\newcommand{\rr}{\operatorname{Re}}
\newcommand{\im}{\operatorname{Im}}
\newcommand{\jj}[2]{J_{#1}(#2)}
\newcommand{\djj}[2]{J'_{#1}(#2)}
\newcommand{\yy}[2]{Y_{#1}(#2)}
\newcommand{\dyy}[2]{Y'_{#1}(#2)}
\newcommand{\hh}[2]{H^{(1)}_{#1}(#2)}
\newcommand{\dhh}[2]{(H^{(1)}_{#1})'(#2)}

\newcommand{\Dom}{\mathop\mathrm{Dom}\nolimits}

\renewcommand{\Re}{\mathop\mathrm{Re}\nolimits}
\renewcommand{\Im}{\mathop\mathrm{Im}\nolimits}

\newcommand{\bigO}[1]{\mathcal{O}{\left(#1\right)}} 

\DeclareMathOperator{\arcsinh}{arcsinh}

\definecolor{DarkGreen}{rgb}{0,0.5,0.1} 

\definecolor{Purple}{rgb}{0.6,0,1}

\makeatother

\begin{document}

\title[]{On Lieb--Thirring inequalities for multidimensional Schr\" odinger operators with complex potentials}

\author{Sabine B{\"o}gli}
\address{
   Department of Mathematical Sciences, Durham University, Upper Mountjoy, Stockton Road, Durham DH1 3LE, UK
    }
\email{sabine.boegli@durham.ac.uk}

\author{Sukrid Petpradittha}
\address{
   Department of Mathematical Sciences, Durham University, Upper Mountjoy, Stockton Road, Durham DH1 3LE, UK
    }
\email{sukrid.petpradittha@durham.ac.uk}

\author{Franti\v sek \v Stampach}
\address{
	Department of Mathematics, Faculty of Nuclear Sciences and Physical Engineering, Czech Technical University in Prague, Trojanova~13, 12000 Praha~2, Czech Republic
	}	
\email{stampfra@cvut.cz}





\begin{abstract}
We solve the open problem by Demuth, Hansmann, and Katriel announced in [Integr. Equ. Oper. Theory 75 (2013), 1--5] by a counter-example construction. The problem concerns a possible generalisation of the Lieb--Thirring inequality for Schr{\" o}dinger operators in $\mathbb{R}^{d}$ to the case of complex-valued potentials. A counter-example has already been found for the one-dimensional case $d=1$ by the first and third authors in [J.~Spectr. Theory 11 (2021), 1391--1413]. Here we generalise the counter-example to higher dimensions $d\geq2$. 
\end{abstract}

\maketitle

\section{Introduction}
Let $p$ depend on the dimension $d$ as follows:
\begin{equation}\label{eq:LiebThCond}
		p\geq 1,\;\mbox{ if }d=1; \qquad 
		p>1,\;\mbox{ if }d=2;		 \qquad
		p\geq d/2,\;\mbox{ if }d\geq3.
\end{equation}
For a real-valued potential $V\in L^p(\mathbb{R}^d)$,  the Schr\" odinger operator $H=-\Delta+V$   is a selfadjoint operator in $L^2(\mathbb{R}^d)$ and its spectrum $\sigma(H)$ is a subset of $\mathbb{R}$. Moreover, $\sigma(H)$ consists of the essential spectrum $\sigma_{\ess}(H)=[0,\infty)$ and the at most countable discrete spectrum $\sigma_{\dd}(H)\subset(-\infty,0)$. The classical Lieb--Thirring inequality states that there exists a constant $C_{p,d}>0$ depending only on $p$ and $d$ such that 
\begin{equation}\label{eq:LiebTh}
	\sum_{\lambda\in\sigma_{\dd}(H)}|\lambda|^{p-d/2}\leq C_{p,d}\|V\|^p_{L^p},
\end{equation}
where the eigenvalues $\lambda$ are counted repeatedly according to their algebraic multiplicities.
Recently, there have been studies on  Lieb--Thirring type inequalities for Schr{\" o}dinger operators $H=-\Delta+V$, where the potential $V\in L^p(\mathbb R^d)$ is allowed to take complex values. For such operators $H$, we still have $\sigma_{\ess}(H)=[0,\infty)$ and $\sigma_{\dd}(H)$ is a set of at most countable isolated eigenvalues of $H$, but these may be non-real.

It turns out that \eqref{eq:LiebTh} does not hold for general complex-valued $V$ in $L^p(\mathbb{R}^d)$ with $p>(d+1)/2$ since, in this case, any point in $\sigma_{\ess}(H)=[0,\infty)$ can be an accumulation point of $\sigma_{\dd}(H)$, see \cite{bogli2017schrodinger,bogli2022counterexample}. 
A possible weaker candidate for a Lieb--Thirring type inequality can be obtained by replacing $|\lambda|^p$ in \eqref{eq:LiebTh} by $\text{dist}(\lambda,[0,\infty))^p$. The resulting inequality then reads
\begin{equation}\label{eq:newLiebTh}
	\sum_{\lambda\in\sigma_{\dd}(H)}\dfrac{\text{dist}(\lambda,[0,\infty))^p}{|\lambda|^{d/2}}\leq C_{p,d}\|V\|^p_{L^p}.
\end{equation}
Note that this inequality reduces to \eqref{eq:LiebTh} when $V$ is real-valued. As it seems to be a reasonable candidate for the Lieb-Thirring inequality extended to complex-valued potentials, the following open question was published in~\cite{demuth2013lieb}.

\vskip8pt

\noindent\textbf{Open Question} (Demuth--Hansmann--Katriel)\textbf{.} \textit{Assuming \eqref{eq:LiebThCond}, does the inequality \eqref{eq:newLiebTh} hold for all $V\in L^p(\mathbb{R}^d)$? Prove or disprove it.}
\medskip

In \cite{bogli2021lieb}, the first and third authors constructed a counter-example in dimension $d=1$ by considering $V$ to be scalar multiples of the characteristic function of the closed interval $[-1,1]$. In \cite[Sec.~3.3]{bogli2021lieb}, we also suggested that the one-parameter family of Schr{\" o}dinger operators
\begin{equation}\label{eq:defH}
H_h:=-\Delta+V_h, \quad V_h:=\ii h\chi_{B_1(0)},
\end{equation}
where $\chi_{B_1(0)}$ is the characteristic function of the unit ball $B_1(0)$ and $h>0$, is a natural candidate for a~counter-example in higher dimensions $d\geq2$. In this article, we show that this is indeed the case. Although it was expected, the extension to the multi-dimensional case is by no means trivial as an involved asymptotic analysis of Bessel functions with complex arguments is needed and the related spectral analysis is in general less explicit than in the one-dimensional case.

Our main result is the following inequality.

\begin{theorem}\label{ThmMainthm}
Let $d\geq 2$, $p>0$, and $0<\varepsilon<1$.
 Then there exists $C_{p,d}>0$ such that, for all sufficiently large $h>0$, we have
\begin{equation}
\frac{1}{\|V_h\|^p_{L^p}}\,\sum_{\lambda\in\sigma_{\dd}(H_h)}\dfrac{{\rm dist}(\lambda,[0,\infty))^p}{|\lambda|^{d/2}}
\geq C_{p,d}\, (\log h)^{\varepsilon}.
\label{eq:main-ineq}
\end{equation}
\end{theorem}

The logarithmic divergence in the parameter~$h$ on the right-hand side of~\eqref{eq:main-ineq} clearly 
answers the question by Demuth, Hansmann and Katriel to the negative.

Recently, the first author proved in \cite{boegli2021improved} a Lieb--Thirring type inequality for Schr{\" o}d\-inger operators with complex-valued potentials $V\in L^p(\mathbb R^d)$ and $p\geq d/2+1$. 
To compensate for the logarithmic divergence, an extra term appears on the left-hand side of~\eqref{eq:newLiebTh} given by a~function of ${\rm dist}(\lambda,[0,\infty))/|\lambda|$; see \cite[Theorem 2.1]{boegli2021improved} for more details. 
This is a generalisation of an earlier result by Demuth, Hansmann and Katriel \cite[Corollary~3]{dem-han-kat_09} which says that for any $0<\tau<1$ there exists $C_{d,p,\tau}>0$ such that 
\begin{equation*}
	\sum_{\lambda\in\sigma_{\dd}(H)}\dfrac{\text{dist}(\lambda,[0,\infty))^{p+\tau}}{|\lambda|^{d/2+\tau}}\leq C_{p,d,\tau}\|V\|^p_{L^p}.
\end{equation*}
Our main result proves that the equation no longer holds if we fomally set $\tau =0$.

A discrete version of Lieb--Thirring type inequalities for Jacobi matrices with complex entries, in particular for one-dimensional discrete Schr{\" o}dinger operators with complex potentials, was found in~\cite{han-kat-guy_11}, and similar open problems published therein have been also answered in~\cite{bogli2021lieb}. Much more literature devoted to Lieb--Thirring inequalities exists nowadays. From those works, whose main focus is on Lieb--Thirring type inequalities for non-self-adjoint Schr{\" o}dinger, Jacobi, and other operators, we mention at least~\cite{bor-gol-kup_09, dem-han-kat_09, dem-han-guy_13, dub_14, fra_18, fra-etal_06, gol-kup_17, GolStep, han_11, sam_14}.

The eigenvalues of our operator $H_h$ are characterised by solutions of a characteristic equation involving Bessel functions. In Section~\ref{sec:char_eq}, we first recall general facts on Schr{\" o}dinger operators with spherically symmetric potentials and then derive the characteristic equation of $H_{h}$. In order to estimate the location and asymptotic behaviour of certain solutions of the characteristic equation, we need to deduce preliminary results concerning the involved Bessel functions, which is worked out in Section~\ref{sec:Bessel}. Finally, in Section~\ref{sec:OnTheProb}, we prove Theorem~\ref{ThmMainthm}.

\section{The characteristic equation}\label{sec:char_eq}

The main goal of this section is to deduce a characteristic equation for the Schr{\" o}dinger operator $H_h$ defined in~\eqref{eq:defH} whose solutions are in direct correspondence with discrete eigenvalues of $H_h$. To do so, we first reduce the eigenvalue equation for $H_h$ to a one-dimensional radial problem. For this step, we recall several well known facts from harmonic analysis of Schr{\" o}dinger operators with spherically symmetric potentials.

\subsection{Schr{\" o}dinger operators with spherically symmetric potentials} 

We consider Schr{\" o}dinger operators in $\R^{d}$ with potentials $V=V(|x|)$, where $|x|$ stands for the Euclidean norm of $x\in\R^{d}$. For the following facts, the reader may consult e.g.\ the book~\cite{stein1971introduction}. Using spherical coordinates in $\R^{d}$, the Laplace operator takes the form
\[
	\Delta=\pd{^2}{r^2}+\dfrac{d-1}{r}\pd{}{r}+\dfrac{1}{r^2}\Delta_{S^{d-1}},
\]
where $r\equiv|x|$ is the radial coordinate and $\Delta_{S^{d-1}}$ is the Laplace--Beltrami operator on the $d$-dimensional unit sphere $S^{d-1}$. The spectrum of $\Delta_{S^{d-1}}$ is discrete, the complete set of eigenfunctions consists of the spherical harmonics $Y^{(\ell)}$ of degree $\ell\in\mathbb{N}_0$, and for each $\ell\in\N_{0}$, we have the eigenvalue equation
\[
	\Delta_{S^{d-1}}Y^{(\ell)}=-\ell(\ell+d-2)Y^{(\ell)}.
\]
The algebraic multiplicity of the eigenvalue $-\ell(\ell+d-2)$ equals
\begin{equation}\label{eq:eigenspace}
	\binom{d+\ell-1}{d-1}-\binom{d+\ell-3}{d-1}
\end{equation}
for each $\ell\in\mathbb{N}_0$, see \cite[p.~140]{stein1971introduction}.

Since the complex-valued potential $V$ depends only on the radius $r\equiv|x|$, we will construct eigenfunctions $\phi\in L^2(\mathbb{R}^d)$ of the Schr\" odinger operator $H=-\Delta+V$ of the form 
\begin{equation}\label{eq:defphi}
	\phi(x)=\psi(r)Y^{(\ell)}\left(\frac{x}{r}\right),
\end{equation}
for some $\ell\in\mathbb{N}_0$ and a function $\psi$ satisfying the radial eigenvalue equation
\begin{equation}\label{eq:eigenvalueProbPsi}
	-\psi''(r)-\dfrac{d-1}{r}\psi'(r)+\dfrac{\ell(\ell+d-2)}{r^2}\psi(r)+V(r)\psi(r)=\lambda\psi(r),
\end{equation}
where $\lambda$ is an eigenvalue of $H$.

\subsection{The characteristic equation}

The potential $V_{h}$ defined in~\eqref{eq:defH} of the Schr{\" o}dinger operator $H_h$ is, of course, spherically symmetric. In the radial coordinate, we have
\[
	V_h(r):=
	\begin{cases}
		\ii h,\quad&\textrm{if $r< 1$},\\
		0,\quad&\textrm{if $r\geq 1$},
	\end{cases}
\]
where $h>0$. By the facts on Schr{\" o}dinger operators with spherically symmetric potentials from the previous subsection, the eigenvalue problem for $H_{h}$ reduces to an analysis of solutions of the one-dimensional eigenvalue equation~\eqref{eq:eigenvalueProbPsi}. Taking into account the special form of the potential $V_{h}$, we seek non-trivial solutions $\psi$ of~\eqref{eq:eigenvalueProbPsi} separately on $(0,1)$ and on $[1,\infty)$ so that $\psi$ and $\psi'$ are continuous at $r=1$ and  
$\psi\in L^2\left((0,\infty),r^{d-1}\dd r\right)$. Then $\phi$ given by~\eqref{eq:defphi} is twice weakly differentiable with $-\Delta\phi=\lambda\phi-V_h\phi \in L^2(\R^d)$, hence $\phi$ belongs to the operator domain $\Dom H_{h}\equiv\Dom(-\Delta)=W^{2,2}(\R^d)$ and so will be an eigenfunction corresponding to the eigenvalue~$\lambda$.

For $0<r<1$, we write $\lambda=k^2$, where $k\in\C$, and introduce a new complex parameter $m\in\C$ such that
\begin{equation}\label{eq:settingLambda}
	\ii h=k^2-m^2.
\end{equation}
Then equation~\eqref{eq:eigenvalueProbPsi} turns into
\begin{align}
-\psi''(r)-\dfrac{d-1}{r}\psi'(r)+\dfrac{\ell(\ell+d-2)}{r^2}\psi(r)-m^2\psi(r)=0,\label{eq:rearrangingPsi}
\end{align}
which has the form of the Bessel differential equation. 
Recall that Bessel's differential equation of order $\nu\in\mathbb{C}$ is the second-order ordinary differential equation
\[
	z^2\dfrac{d^2w}{dz^2}+z\dfrac{dw}{dz}+(z^2-\nu^2)w=0,
\]
and has solutions called the Bessel functions of the first kind $\jj{\pm\nu}{z}$, of the second kind $\yy{\nu}{z}$,  and of the third kind $\hh{\nu}{z},H^{(2)}_\nu(z)$ (also called Hankel functions). Of course, these solution are interrelated. As the main source for the theory of Bessel functions, we use the classical book by Watson~\cite{watson1995treatise}. For a more updated and well arranged source of numerous formulas for Bessel functions we use the digital library~\cite{dlmf} which replaced the older book~\cite{abramowitz1948handbook} of Abramowitz and Stegun. Let us also mention Olver's book~\cite{olver1997asymptotics}, where proofs on various asymptotic formulas for Bessel functions can be found.

For $0<r<1$, this leads us to consider the solution $\psi$ of~\eqref{eq:rearrangingPsi} having the form 
\[
	\psi(r)=r^{1-d/2}J_{\nu}(mr),
\]
where
\begin{equation}\label{eq:Nudefinition}
	\nu=\ell+\frac{d}{2}-1.
\end{equation}
Since for $\nu\geq0$, the function $\jj{\nu}{z}$ remains bounded as $z\to 0$ in the half-plane $\Re z>0$, see e.g.~\cite[Eq.~(10.7.3)]{dlmf}, we have 
\[
\psi\in L^2\left((0,1), r^{d-1}\dd r\right).
\] 
Later on, see Remark~\ref{rem:fourth-quadrant}, our analysis will be confined to the fourth quadrant of $\C$ in the variable $m$, i.e. $\Re m>0$ and $\Im m<0$, hence we may assume $\Re z>0$ above.

For $r\geq 1$, $V(r)=0$, and by writing $\lambda=k^2\in\mathbb{C}$ in equation \eqref{eq:eigenvalueProbPsi}, we obtain
\[
	-\psi''(r)-\dfrac{d-1}{r}\psi'(r)+\dfrac{\ell(\ell+d-2)}{r^2}\psi(r)-k^2\psi(r)=0.
\]
To ensure that 
\[
\psi\in L^2\left((1,\infty),r^{d-1}\dd r\right),
\]
we choose the Hankel function of the first kind $H^{(1)}_\nu(z)$ as a solution in the corresponding Bessel equation because $H^{(1)}_\nu(z)$ is exponentially decaying as $z\to\infty$ provided that $\Im z>0$, see \cite[Eq.~(10.2.5)]{dlmf}. Consequently, up to a multiplicative constant, we may take
\[
	\psi(r)=r^{1-d/2}H^{(1)}_\nu(kr)
\]
for $r>1$ provided that 
\begin{equation}
	\im k>0.
\label{eq:im_k_pos}
\end{equation}

Now, combing both solutions on $(0,1)$ and $[1,\infty)$, and choosing multiplicative constants so that the resulting function is continuous at $r=1$, we get the solution 
\[
	\psi(r)=	
	\begin{cases}
		H^{(1)}_\nu(k)r^{1-d/2}J_{\nu}(mr),\quad&\mbox{ if } 0<r<1,\\
		J_\nu(m)r^{1-d/2}H^{(1)}_\nu(kr),\quad&\mbox{ if } r\geq 1,
	\end{cases}
\]
provided that we also have~\eqref{eq:im_k_pos}.

Lastly, in order for the constructed solution $\psi$ to become the radial eigenfunction of $H_h$, its first derivative $\psi'$ must be continuous at $r=1$, too. It follows that a last condition needs to be imposed:
\[
	mJ_\nu'(m)H^{(1)}_\nu(k)-kJ_\nu(m)(H^{(1)}_\nu)'(k)=0.
\]
It means that
\begin{equation}\label{eq:CharacteristicSign}
	\frac{k}{m}=\dfrac{\djj{\nu}{m}\hh{\nu}{k}}{\jj{\nu}{m}\dhh{\nu}{k}},
\end{equation}
if no division by $0$ occurs.
In the proof below, we will restrict our parameter region to complex numbers for which we can exclude the occurrence of division by $0$. 
Equation~\eqref{eq:CharacteristicSign} is the desired characteristic equation.
In the following it turns out to be easier to work with the squared version of the latter equation. 
Using \eqref{eq:settingLambda}, it can be written as
\begin{equation}\label{eq:Characteristic}
	\frac{\ii h}{m^{2}}+1=\left(\dfrac{\djj{\nu}{m}\hh{\nu}{k}}{\jj{\nu}{m}\dhh{\nu}{k}}\right)^{\!2}.
\end{equation}
Now, for each $h>0$ and $\ell\in\mathbb{N}_0$ (determining $\nu$ by~\eqref{eq:Nudefinition}), if any 
$m, k\in\C$ with $\Re m>0$ and $\Im k>0$ 
satisfy~\eqref{eq:CharacteristicSign},
then $\lambda=k^{2}$ is a discrete eigenvalue of $H_h$ of algebraic multiplicity at least~\eqref{eq:eigenspace}.
In the proof below, we first find solutions of \eqref{eq:Characteristic} and then check that these solutions in fact solve~\eqref{eq:CharacteristicSign}, i.e.\ that the signs on both sides of~\eqref{eq:CharacteristicSign} agree.

We conclude this section by locating eigenvalues of $H_{h}$ in a strip, which will be used later.

\begin{lemma}\label{LemmaLocationEigenvalues}
	Let $h>0$ and $\lambda$ be an eigenvalue of $H_h$. Then $\rr\lambda\geq 0$ and $0\leq\im\lambda\leq h$.
\end{lemma}
\begin{proof}
We show that the numerical range $W(H_{h})$ of $H_{h}$ is a subset of the strip determined by the equations $\rr\lambda\geq 0$ and $0\leq\im\lambda\leq h$. Since the point spectrum of $H_{h}$ is a subset of $W(H_{h})$ the claim follows.

Suppose $\phi\in\Dom H_h$ to be normalised, $\|\phi\|_{L^2}=1$. Clearly,
	\[
\left\langle{H_h\phi,\phi}\right\rangle=\left\langle{-\Delta \phi,\phi}\right\rangle+\ii h\left\langle{\chi_{B_1(0)}\phi,\phi}\right\rangle.
	\]
As the Laplacian is non-negative, we have
	\[
		\left\langle{-\Delta \phi,\phi}\right\rangle 
		=\int_{\mathbb{R}^d}|\nabla \phi(x)|^2\dd x\geq 0.
	\]
Moreover, 
	\[
		\left\langle{\chi_{B_1(0)}\phi,\phi}\right\rangle=\int_{B_1(0)}|\phi(x)|^2\,\dd x\leq\int_{\mathbb{R}^d}|\phi(x)|^2\,\dd x=1,
	\]
by assumption. It follows that $W(H_{h})\subset[0,\infty)+\ii[0,h]$ for $h>0$, and the proof is complete.
\end{proof}

\section{Selected properties of Bessel functions}\label{sec:Bessel}

In order to analyse zeros of the  equation~\eqref{eq:Characteristic} for $h$ large, we need to understand the asymptotic behaviour of the involved Bessel functions in a special regime when both the main argument as well as the order tend simultaneously but not independently to infinity. 
This section collects auxiliary results on selected properties of Bessel functions that will be used later in the proof of Theorem~\ref{ThmMainthm}.

\subsection{Two expansion formulas for Bessel functions}

For $n\in\mathbb{N}_{0}$, we adopt Hankel's notation:
\[
(\nu,n):=\frac{1}{n!}\prod_{j=1}^{n}\left[\nu^{2}-\left(j-\frac{1}{2}\right)^{\!2}\right],
\]
see~\cite[p.~198]{watson1995treatise}. Our proof heavily relies on the following expansions of Bessel functions with a certain uniform control of remainder terms in a complex sector.

\begin{lemma}\label{lem:HankelAsymp}
	Let $p\in\N_0$. Then for all $\nu\geq p$ and $z\in\C$ satisfying $0<|z|\leq 2\Re z$, we have 
\begin{equation}
\jj\nu{z}\djj\nu{z}+\yy\nu{z}\dyy\nu{z}=-\dfrac{1}{\pi}\sum_{n=0}^{p-1}\frac{(2n+1)!}{4^{n}n!}\frac{(\nu,n)}{z^{2n+2}}+\mathcal{X}_{p}(z;\nu)
\label{eq:HankelAsymp1}
\end{equation}
and 
\begin{equation}
		J^2_\nu(z)+Y^2_\nu(z)=\dfrac{2}{\pi}\sum_{n=0}^{p-1}\frac{(2n)!}{4^{n}n!}\frac{(\nu,n)}{z^{2n+1}}+\mathcal{Y}_{p}(z;\nu)
\label{eq:HankelAsymp2}
\end{equation}
with
\[
 |\mathcal{X}_{p}(z;\nu)|\leq C_{p}\,\frac{\nu^{2p}}{| z|^{2p+2}}
 \quad\mbox{ and }\quad
  |\mathcal{Y}_{p}(z;\nu)|\leq C_{p}'\,\frac{\nu^{2p}}{| z|^{2p+1}},
\]
where the constants $C_{p}, C_{p}'>0$ depend on $p$, but are independent of $\nu$ and $z$.
\end{lemma}

\begin{proof}
We prove~\eqref{eq:HankelAsymp1} together with the inequality for the remainder term. The proof of~\eqref{eq:HankelAsymp2} is similar and will only be indicated at the end.
 
Our starting point is Nicholson's integral representation formula
\begin{equation}\label{eq:Nicholson}
J^2_\nu(z)+Y^2_\nu(z)=\dfrac{8}{\pi^2}\int_0^\infty K_0(2z\sinh(t))\cosh (2\nu t)\,\dd t,
\end{equation}
which holds true if $\rr{z}>0$, see \cite[Eq.~(10.9.30)]{dlmf} or \cite[p. 444]{watson1995treatise}. The function $K_0$ in the integrand is the modified Bessel function of the second kind of order zero. By differentiating~\eqref{eq:Nicholson} with respect to $z$ and using the identity $K_0'(z)=-K_1(z)$ valid for all $z\in\C$, see \cite[Eq.~(10.29.3)]{dlmf}, one obtains the formula
\begin{equation}\label{eq:diff_J^2+Y^2WithIntegral}
	\jj\nu{z}\djj\nu{z}+\yy\nu{z}\dyy\nu{z}=-\dfrac{8}{\pi^2}\int_0^\infty K_1(2z\sinh(t))\cosh (2\nu t)\sinh(t)\,\dd t.
\end{equation}

Next, we make use of the expansion
\begin{equation}\label{eq:errorExpansion}
	\dfrac{\cosh(2\nu t)}{\cosh(t)}=\sum_{n=0}^{p-1}\dfrac{n!(\nu,n)}{(2n)!}2^{2n}\sinh^{2n}(t)+R_p(t),
\end{equation}
where
\[
|R_p(t)|\leq\left|\dfrac{\cos(\nu\pi)}{\cos(\rr{(\nu\pi)})}\right|\dfrac{p!|(\rr{\nu},p)|}{(2p)!}2^{2p}\sinh^{2p}(t),
\]
see \cite[p. 448]{watson1995treatise}. As our $\nu$ is real, the upper bound for the remainder term $R_{p}$ simplifies to
\begin{equation}\label{eq:errorBound}
	|R_p(t)|\leq \dfrac{p!|(\nu,p)|}{(2p)!}2^{2p}\sinh^{2p}(t).
\end{equation}
Plugging \eqref{eq:errorExpansion} into \eqref{eq:diff_J^2+Y^2WithIntegral} yields
\begin{equation}
	\jj\nu{z}\djj\nu{z}+\yy\nu{z}\dyy\nu{z}=-\dfrac{8}{\pi^2}\left(I_{1,p}+I_{2,p}\right),
\label{eq:wronsk_I1_I2_aux}
\end{equation}
where
\[
I_{1,p}:=\sum_{n=0}^{p-1}\dfrac{n!(\nu,n)}{(2n)!}2^{2n}\int_0^\infty K_1(2z\sinh(t))\sinh^{2n+1}(t)\cosh(t)\,\dd t,
\]
and
\[
I_{2,p}:=\int_0^\infty K_1(2z\sinh(t))\sinh(t)\cosh(t)R_p(t)\,\dd t.
\]

Further we consider separately the integrals $I_{1,p}$ and $I_{2,p}$, starting with $I_{1,p}$. By the change of variable $u=\sinh t$, we obtain
\[
I_{1,p}=\sum_{n=0}^{p-1}\dfrac{n!(\nu,n)}{(2n)!}2^{2n}\int_0^\infty K_1(2zu)u^{2n+1}\, \dd u.
\]
The last integral can be changed to a complex contour integral along a ray from $0$ to complex infinity of angle $\arg z$ that is located in the right half-plane since $\Re z>0$. 
By the analyticity of $K_{1}$ and the standard homotopy argument, taking also into account the asymptotic behaviour of $K_{1}(z)$ as $z\to0$ and $z\to\infty$ from the right half-plane $\Re z>0$, see~\cite[Eqs.~(10.30.2) and~(10.25.3)]{dlmf}, we may deform the integration contour to the positive half-line getting the equality
\[
\int_0^\infty K_1(2zu)u^{2n+1}\, \dd u=\frac{1}{(2z)^{2n+2}}\int_0^{\infty} K_1(t)t^{2n+1}\, \dd t.
\]
Therefore, we have
\[
I_{1,p}=\frac{1}{4}\sum_{n=0}^{p-1}\dfrac{n!}{(2n)!}\frac{(\nu,n)}{z^{2n+2}}\int_0^\infty K_1(t)t^{2n+1}\, \dd t.
\]
Next we apply the integral identity 
\begin{equation}\label{eq:Kint}
\int_0^\infty K_{\mu_1}(t)t^{\mu_2-1}\dd t=2^{\mu_2-2}\Gamma\left(\dfrac{\mu_2-\mu_1}{2}\right)\Gamma\left(\dfrac{\mu_2+\mu_1}{2}\right),
\end{equation}
which is valid when $\rr\mu_2>|\rr\mu_1|$, see \cite[p. 388, Eq. (8)]{watson1995treatise}, and arrive at the expression
\begin{equation}
I_{1,p}=\frac{1}{4}\sum_{n=0}^{p-1}\frac{2^{2n} n!}{(2n)!}\frac{(\nu,n)}{z^{2n+2}}\Gamma\left(n+\frac{1}{2}\right)\Gamma\left(n+1+\dfrac{1}{2}\right)=\frac{\pi}{8}\sum_{n=0}^{p-1}\frac{(2n+1)!}{2^{2n}n!}\frac{(\nu,n)}{z^{2n+2}},
\label{eq:I1_fin_expr_aux}
\end{equation}
after a simplification of the Gamma functions.

As a next step, we estimate $|I_{2,p}|$ from above. Comparing~\eqref{eq:wronsk_I1_I2_aux} and~\eqref{eq:I1_fin_expr_aux} with~\eqref{eq:HankelAsymp1}, we already see that $\mathcal{X}_{p}(z;\nu)=-8I_{2,p}/\pi^{2}$. Again, the change of variable $u=\sinh t$ yields
\[
 I_{2,p}=\int_0^\infty u K_1(2zu)R_p(\arcsinh u)\, \dd u.
\]
Then, by~\eqref{eq:errorBound}, we get
\begin{equation}
|I_{2,p}|\leq\int_0^\infty u|K_1(2zu)||R_p(\arcsinh u)|\,\dd u\leq \frac{2^{2p}p!}{(2p)!}\,|(\nu,p)|\int_0^\infty |K_1(2zu)|u^{2p+1}\,\dd u.
\label{eq:mod_I_2p_aux}
\end{equation}
In order to deduce a suitable upper bound for $|K_1(2zu)|$, we employ the integral representation
\[
K_{1}(z)=\int_0^\infty \e^{-z\cosh(x)}\cosh(x)\, \dd x
\]
which is valid if $\Re z>0$, see~\cite[Eq.~(10.32.9)]{dlmf}. At this point, we use the assumption $|z|\leq 2\Re z$ and deduce the estimate
\[
|K_1(2zu)|\leq\int_0^\infty \e^{-2(\rr{z})u\cosh(x)}\cosh(x)\,\dd x\leq\int_0^\infty \e^{-|z|u\cosh(x)}\cosh(x)\,\dd x=K_1(u|z|).
\]
By using the last estimate in~\eqref{eq:mod_I_2p_aux}, changing the variable $t=u|z|$, applying identity~\eqref{eq:Kint}, and simplifying the resulting expression with Gamma functions similarly as above, one arrives at the estimate
\[
 |I_{2,p}|\leq \frac{\pi}{4}\frac{(2p+2)!}{(p+1)!}\frac{|(\nu,p)|}{|z|^{2p+2}}.
\]
Lastly, as we assume that $\nu\geq p$, we may also estimate the Hankel bracket
\[
 |(\nu,p)|=\frac{1}{p!}\prod_{j=1}^{p}\left[\nu^{2}-\left(j-\frac{1}{2}\right)^{\!2}\right]\leq\frac{\nu^{2p}}{p!}.
\]
In total, we have
\[
|I_{2,p}|\leq \frac{\pi}{4}\frac{(2p+2)!}{p!(p+1)!}\frac{\nu^{2p}}{|z|^{2p+2}}.
\]
Putting together equations~\eqref{eq:wronsk_I1_I2_aux} and \eqref{eq:I1_fin_expr_aux} with the last estimate, we prove~\eqref{eq:HankelAsymp1} together with the inequality for the remainder term $\mathcal{X}_{p}(z;\nu)$.

The proof of the second formula~\eqref{eq:HankelAsymp2} proceeds similarly. Instead of~\eqref{eq:diff_J^2+Y^2WithIntegral}, we expand the term $\cosh(2\nu t)$ using~\eqref{eq:errorExpansion} directly in Nicholson's integral~\eqref{eq:Nicholson}. Then we continue analogously as before to simplify the integral as in case of $I_{1,p}$ and estimate the error term as in case of $I_{2,p}$ getting the identity~\eqref{eq:HankelAsymp2} as well as the upper bound for $|\mathcal{Y}_{p}(z;\nu)|$. The proof of Lemma~\ref{lem:HankelAsymp} is complete.
\end{proof}

\subsection{The modulus and phase functions}

Further elements of the theory of Bessel functions, which will be frequently used, involve the so-called modulus and phase functions $M_{\nu}$ and $\theta_{\nu}$ that are introduced by means of the equations 
\begin{equation}
 J_{\nu}(z)=\sqrt{\frac{2}{\pi z}}\,M_{\nu}(z)\cos\theta_{\nu}(z)
 \quad\mbox{ and }\quad 
 Y_{\nu}(z)=\sqrt{\frac{2}{\pi z}}\,M_{\nu}(z)\sin\theta_{\nu}(z).
\label{eq:J-Y_M-psi}
\end{equation}
In general, the equations in~\eqref{eq:J-Y_M-psi} give rise to multi-valued functions $M_{\nu}$ and $\theta_{\nu}$. The standard branch of $M_{\nu}$ is determined by requiring $M_{\nu}$ to be continuous in $(0,\infty)$ and $M_{\nu}(x)>0$ for $x>0$, and as such, $M_{\nu}$ extends to an analytic function in the right half-plane $\Re z>0$. It follows immediately from~\eqref{eq:J-Y_M-psi} that 
\begin{equation}
 M_{\nu}^{2}(z)=\frac{\pi z}{2}\left(J_{\nu}^{2}(z)+Y_{\nu}^{2}(z)\right).
 \label{eq:M_id}
\end{equation}
From the form of the equations in~\eqref{eq:J-Y_M-psi}, one may suggest to incorporate the factor $\sqrt{2/(\pi z)}$ into the modulus function $M_{\nu}(z)$. This alternative notation for the modulus function, which is unfortunately denoted by $M_{\nu}$ again, is also used in more modern literature on Bessel functions, see for example~\cite[\S~10.18]{dlmf}. Here we stick with the original notation of Marshall for $M_{\nu}$ from~\eqref{eq:J-Y_M-psi} that is also used in Watson's treatise~\cite{watson1995treatise}. One needs to be careful with this double meaning of $M_{\nu}$ when, for example, formulas listed in~\cite[\S~10.18]{dlmf} are used; cf.~\eqref{eq:M_id} with~\cite[Eq.~(10.18.6)]{dlmf}.

It also readily follows from~\eqref{eq:J-Y_M-psi} that 
\begin{equation}
  \theta_{\nu}(z)=\arctan\frac{Y_{\nu}(z)}{J_{\nu}(z)},
 \label{eq:phase_id}
\end{equation}
cf.~\cite[Eq.~(10.18.7)]{dlmf}, from which, using the Wronskian identity~\cite[p.~76]{watson1995treatise}
\begin{equation}
 \jj\nu{z}\dyy\nu{z}-\djj\nu{z}\yy\nu{z}=\frac{2}{\pi z},
\label{eq:wronsk_JY}
\end{equation}
one infers that
\begin{equation}
  \theta_{\nu}'(z)=\frac{J_{\nu}(z)Y_{\nu}'(z)-J_{\nu}'(z)Y_{\nu}(z)}{J_{\nu}^{2}(z)+Y_{\nu}^{2}(z)}=\frac{2}{\pi z}\frac{1}{{J_{\nu}^{2}(z)+Y_{\nu}^{2}(z)}}.
 \label{eq:phase-der_id}
\end{equation}
The standard branch of $\theta_{\nu}$, which is adopted below, is determined by requiring $\theta_{\nu}$ to be continuous in $(0,\infty)$ and $\theta_{\nu}(x)\to-\pi/2$ as $x\to0+$. Below, we verify, for $\nu\geq1$, that this branch of $\theta_{\nu}$ is an analytic function in the open convex set
\begin{equation}
 \mathcal{M}_{\nu}:=\{z\in\C \mid  A\nu<\Re z \mbox{ and } |z|<2\Re z\},
 \label{eq:set_M}
\end{equation}
with 
\begin{equation}\label{eq:defA}
A:=\sqrt{\pi C_{1}'},
\end{equation}
where $C_{1}'>0$ is the constant from Lemma~\ref{lem:HankelAsymp} with $p=1$. This is, in fact, a subset of a half-plane where the standard branch of $\theta_{\nu}$ extends to a single-valued   analytic  function, but this is not needed for our purposes. In addition to the analyticity region for the phase function, the definition of the set $\mathcal{M}_{\nu}$ is also designed for other aspects of the forthcoming analysis.

The modulus and phase functions were introduced by Marshall and used in an asymptotic analysis of large zeros of the cylinder function, see \cite[p. 505]{watson1995treatise}. In fact, Marshall introduced a~slightly different phase function given by 
\begin{equation}
 \psi_{\nu}(z):=z-\frac{\nu\pi}{2}-\frac{\pi}{4}-\theta_{\nu}(z),
\label{eq:def_psi}
\end{equation}
and worked with a different pair of functions $P(\cdot,\nu)$ and $Q(\cdot,\nu)$ related to $J_{\nu}$ and $Y_{\nu}$; the reader may find details in~\cite{watson1995treatise}. We also occasionally use the function $\psi_{\nu}$ below. The modulus and phase functions were used in more recent work~\cite{heitman2015asymptotics} to study the asymptotic behaviour of $J_{\nu}(z)$ and $Y_{\nu}(z)$ in the so-called Fresnel regime $|z|>\nu$, but with $z$ confined to positive reals. This is close but not exactly what is needed in our analysis. Although not directly applied, the methods from paper~\cite{heitman2015asymptotics} as well as Marshall's original ideas inspired our approach.

\begin{lemma}\label{lem:phase_analytic}
 For $\nu\geq1$, the phase functions $\theta_{\nu}$ and $\psi_{\nu}$ are analytic in the set $\mathcal{M}_{\nu}$ defined by~\eqref{eq:set_M}. Moreover, we have the estimates
\begin{equation}
  |\psi_{\nu}'(z)|\leq A^{2}\frac{\nu^{2}}{|z|^{2}}
 \label{eq:psi_der_bound}
\end{equation}
and
\begin{equation}
|\psi_{\nu}''(z)|\leq\tilde{A}\frac{\nu^{2}}{|z|^{3}}
 \label{eq:psi_der2_bound}
\end{equation}
for all $z\in\mathcal{M}_{\nu}$, where $A$ is the constant in \eqref{eq:defA}, $\tilde{A}:=2\pi(C_{1}'+2C_{1})>0$, and $C_1, C_1'$ are the constants from Lemma~\ref{lem:HankelAsymp}.
\end{lemma}

\begin{proof}
 First we verify the analyticity of $\theta_{\nu}$ in $\mathcal{M}_{\nu}$. Then $\psi_{\nu}$ must be analytic in $\mathcal{M}_{\nu}$, too, by its definition~\eqref{eq:def_psi}. Bearing~\eqref{eq:phase-der_id} in mind, the value of the standard branch of $\theta_{\nu}$ at $z\in\mathcal{M}_{\nu}$ can be expressed by the contour integral
 \[
  \theta_{\nu}(z)=-\frac{\pi}{2}+\frac{2}{\pi}\int_{0}^{z}\frac{1}{{J_{\nu}^{2}(u)+Y_{\nu}^{2}(u)}}\frac{\dd u}{u}
 \]
 where the integration contour consists of the interval $(0,A\nu)$ and then continues arbitrarily to the point $z$ within the set $\mathcal{M}_{\nu}$, provided that the analytic function $J_{\nu}^{2}+Y_{\nu}^{2}$ has no zeros in $\mathcal{M}_{\nu}$. Notice that the integrand is indeed integrable near $0+$ for $\nu>0$, since $J_{\nu}(u)=\bigO{u^{\nu}}$ and $Y_{\nu}(u)=\bigO{u^{-\nu}}$ as $u\to0+$, see~\cite[Eqs.~(10.7.3--4)]{dlmf}.
 
Once we check that $J_{\nu}^{2}(z)+Y_{\nu}^{2}(z)\neq0$ for all $z\in\mathcal{M}_{\nu}$, the analyticity of $\theta_{\nu}$ in $\mathcal{M}_{\nu}$ follows from the contour integral representation. To this end, we apply Lemma~\ref{lem:HankelAsymp} with $p=1$ getting, for any $z\in\mathcal{M}_{\nu}$, the estimate
\[
 |J_{\nu}^{2}(z)+Y_{\nu}^{2}(z)|\geq\frac{2}{\pi|z|}-|\mathcal{Y}_{1}(z;\nu)|\geq\frac{2}{\pi|z|}-C_{1}'\frac{\nu^{2}}{|z|^{3}}=\frac{2}{\pi|z|}-\frac{A^{2}\nu^{2}}{\pi|z|^{3}}\geq\frac{1}{\pi|z|}>0,
\]
where we have used that $A\nu<\Re z\leq|z|$ for $z\in\mathcal{M}_{\nu}$. Consequently, $\mathcal{M}_{\nu}$ is a zero-free region for $J_{\nu}^{2}+Y_{\nu}^{2}$, indeed.

Second we establish~\eqref{eq:psi_der_bound}. Substituting from~\eqref{eq:HankelAsymp2} with $p=1$ in~\eqref{eq:phase-der_id}, we get
\begin{equation}\label{eq:part1}
 \theta_{\nu}'(z)=\frac{2}{\pi z}\frac{1}{{J_{\nu}^{2}(z)+Y_{\nu}^{2}(z)}}=\frac{2}{2+\pi z\mathcal{Y}_{1}(z;\nu)}.
\end{equation}
Using also~\eqref{eq:def_psi}, 
we obtain the expression
\[
\psi_{\nu}'(z)=1-\theta_{\nu}'(z)=\frac{\pi z\mathcal{Y}_{1}(z;\nu)}{2+\pi z\mathcal{Y}_{1}(z;\nu)}.
\]
By Lemma~\ref{lem:HankelAsymp} and definition~\eqref{eq:set_M} of the set $\mathcal{M}_{\nu}$, we have
\begin{equation}\label{eq:part2}
\left|\pi z\mathcal{Y}_{1}(z;\nu)\right|\leq\pi C_{1}'\frac{\nu^{2}}{|z|^{2}}=A^{2}\frac{\nu^{2}}{|z|^{2}}\leq1
\end{equation}
for all $z\in\mathcal{M}_{\nu}$. Thus, for any $z\in\mathcal{M}_{\nu}$, we may estimate as
\[
 |\psi_{\nu}'(z)|\leq\frac{|\pi z\mathcal{Y}_{1}(z;\nu)|}{2-|\pi z\mathcal{Y}_{1}(z;\nu)|}\leq A^{2}\frac{\nu^{2}}{|z|^{2}},
\]
proving~\eqref{eq:psi_der_bound}.

Lastly, we prove~\eqref{eq:psi_der2_bound}. By \eqref{eq:def_psi}, 
we have $\psi_{\nu}''=-\theta_{\nu}''$. Differentiating~\eqref{eq:phase-der_id} once more, we find
\[
 \psi_{\nu}''(z)=\frac{2}{\pi z}\frac{1}{J_{\nu}^{2}(z)+Y_{\nu}^{2}(z)}\left[\frac{1}{z}+\frac{2J_{\nu}(z)J_{\nu}'(z)+2Y_{\nu}(z)Y_{\nu}'(z)}{J_{\nu}^{2}(z)+Y_{\nu}^{2}(z)}\right].
\]
Bearing in mind~\eqref{eq:phase-der_id} and using Lemma~\ref{lem:HankelAsymp} with $p=1$, we arrive at the expression
\[
 \psi_{\nu}''(z)=\theta_{\nu}'(z)\frac{\pi\mathcal{Y}_{1}(z;\nu)+2\pi z\mathcal{X}_{1}(z;\nu)}{2+\pi z\mathcal{Y}_{1}(z;\nu)}.
\]
For any $z\in\mathcal{M}_{\nu}$, we may estimate roughly $|\theta_{\nu}'(z)|\leq 2$ by \eqref{eq:part1} and \eqref{eq:part2}.
The remaining fraction is to be estimated with the aid of Lemma~\ref{lem:HankelAsymp} as follows:
\[
 \left|\frac{\pi\mathcal{Y}_{1}(z;\nu)+2\pi z\mathcal{X}_{1}(z;\nu)}{2+\pi z\mathcal{Y}_{1}(z;\nu)}\right|\leq\frac{\pi|\mathcal{Y}_{1}(z;\nu)|+2\pi|z\mathcal{X}_{1}(z;\nu)|}{2-\pi| z\mathcal{Y}_{1}(z;\nu)|}\leq\pi(C_{1}'+2C_{1})\frac{{\nu}^{2}}{|z|^{3}}.
\]
This completes the proof of inequality~\eqref{eq:psi_der2_bound}.
\end{proof}

From the bounds on the derivative of $\psi_{\nu}$ from Lemma~\ref{lem:phase_analytic}, we further deduce bounds on the function $\psi_{\nu}$ itself in $\mathcal{M}_{\nu}$.

\begin{lemma}\label{lem:PhaseAsymp}
 Let $\nu\geq1$. Then for all $z\in\mathcal{M}_{\nu}$, we have the inequalities
 \begin{equation}\label{eq:PhaseAsymp}
 |\psi_{\nu}(z)|\leq 2A^{2}\frac{\nu^{2}}{|z|}
 \end{equation}
 and
 \begin{equation}\label{eq:imaginaryPhase}
 |\im\psi_\nu(z)|\leq 4A^{2} |\im z|\frac{\nu^2}{|z|^2},
 \end{equation}
 where $A>0$ is the constant in \eqref{eq:defA}. 
\end{lemma}


\begin{proof}
First we verify~\eqref{eq:PhaseAsymp}. The idea of the proof is to integrate the expansion for the derivative $\psi_{\nu}'$ and determine the constant term. First note that if $z\in\mathcal{M}_{\nu}$, then $z+[0,\infty)\subset\mathcal{M}_{\nu}$. We show that, with $\nu\geq1$ fixed,
\begin{equation}
\psi_{\nu}(z)\to0 \mbox{ as } z\to\infty \mbox{ along a horizontal ray in } \mathcal{M}_{\nu}.
\label{eq:psi_lim_zero}
\end{equation} 
For $z\in\mathcal{M}_{\nu}$ it follows that
\[
 |\psi_{\nu}(z)|\leq|\psi_{\nu}(z)-\psi_{\nu}(\Re z)|+|\psi_{\nu}(\Re z)|.
\]
Due to~\cite[Eq.~(10.18.18)]{dlmf} with definition~\eqref{eq:def_psi}, 
$|\psi_{\nu}(\Re z)|\to 0$ as $\Re z\to\infty$. Now, to determine the limit of $|\psi_{\nu}(z)-\psi_{\nu}(\Re z)|$, we make use of the analyticity of $\psi_{\nu}$ in $\mathcal{M}_{\nu}$, see Lemma~\ref{lem:phase_analytic}, and the mean value theorem  to infer that
\[
 |\psi_{\nu}(z)-\psi_{\nu}(\Re z)|\leq|\Im z|\sup_{\xi\in(\Re z,z)}|\psi_{\nu}'(\xi)|,
\]
where $(\Re z,z)$ is the open vertical line segment connecting the points $\Re z$ and $z$. Since this line segment is entirely contained in $\mathcal{M}_{\nu}$, the inequality~\eqref{eq:psi_der_bound} is applied to deduce
\[
 |\psi_{\nu}(z)-\psi_{\nu}(\Re z)|\leq|\Im z|\frac{A^{2}\nu^{2}}{(\Re z)^{2}}\to 0
\]
as $\Re z\to\infty$ with $\Im z$ constant. The claim from~\eqref{eq:psi_lim_zero} follows.

Next, we consider the ray $z+[0,\infty)$ for 
$z\in\mathcal{M}_{\nu}$. As said before, this ray is
entirely located in the set $\mathcal{M}_{\nu}$, on which $\psi_{\nu}$ is analytic. Consequently, bearing also~\eqref{eq:psi_lim_zero} in mind, we may express $\psi_{\nu}(z)$ as the contour integral
\[
 \psi_{\nu}(z)=-\int_{z}^{+\infty+\ii \Im z}\psi_\nu'(u)\,\dd u.
\]
Notice that the integral converges since $|\psi_\nu'(u)|$ decays at least as $1/|u|^{2}$ for $\Re u$ large, as we know from~\eqref{eq:psi_der_bound}. Then, with the aid of~\eqref{eq:psi_der_bound}, we may readily estimate
\[
|\psi_{\nu}(z)|=\int_{z}^{+\infty+\ii \Im z}|\psi_\nu'(u)|\,|\dd u|\leq A^{2}\nu^{2}\int_{\Re z}^{\infty}\frac{\dd x}{x^{2}}= \frac{A^{2}\nu^{2}}{\Re z}\leq 2A^{2}\,\frac{\nu^{2}}{|z|},
\]
proving~\eqref{eq:PhaseAsymp}.

Lastly, we verify~\eqref{eq:imaginaryPhase}. Fix $z\in\mathcal{M}_{\nu}$. Recall that $\psi_{\nu}(x)\in\mathbb{R}$ if $x>0$. Then, by the above arguments, we have 
\[
 |\Im\psi_{\nu}(z)|=|\Im\psi_{\nu}(z)-\Im\psi_{\nu}(\Re z)|\leq|\Im z|\sup_{\xi\in(\Re z,z)}|\psi_{\nu}'(\xi)|
 \leq|\Im z|\frac{A^{2}\nu^{2}}{(\Re z)^{2}}.
\]
We employ $|z|<2\Re z$ to argue
\[
 |\Im\psi_{\nu}(z)|\leq|\Im z|\frac{4A^{2}\nu^{2}}{|z|^{2}}.
\]
The proof is complete.
\end{proof}

We will also prepare formulas for an asymptotic analysis of the logarithmic derivatives of the functions $J_{\nu}$ and $H_{\nu}^{(1)}$ that appear in the  equation~\eqref{eq:Characteristic}.

\begin{lemma}\label{lem:log_der}
 Let $\nu\geq1$ and $z\in\mathcal{M}_{\nu}$. Then we have
 \begin{equation}
  \frac{\djj{\nu}{z}}{\jj{\nu}{z}}=-\theta_{\nu}'(z)\tan\theta_\nu(z)+\mathcal{Z}_{0}(z;\nu),
 \label{eq:log_der_J}
 \end{equation}
 where
 \[
  |\mathcal{Z}_{0}(z;\nu)|\leq\frac{C_{0}\pi}{|z|}.
 \]
 If, in addition, $|z|>1+C_{1}/C_{1}'$, we have
 \begin{equation}
 	\theta_{\nu}'(z)\frac{\hh{\nu}{z}}{\dhh{\nu}{z}}=-\ii+\mathcal{Z}_{1}(z;\nu),
  \label{eq:log_der_H}
 \end{equation}
 where 
 \[
  |\mathcal{Z}_{1}(z;\nu)|\leq\left(1+\frac{C_1}{C_1'}\right)\frac{1}{|z|}.
 \]
 The positive constants $C_{0}$, $C_{1}$, and $C_{1}'$ are those of Lemma~\ref{lem:HankelAsymp}.
\end{lemma}

\begin{proof}
We verify~\eqref{eq:log_der_J}.
First, taking the logarithmic derivative with respect to $z$ in the first equation of~\eqref{eq:J-Y_M-psi} yields the formula
\[
\frac{\djj{\nu}{z}}{\jj{\nu}{z}}=-\frac{1}{2z}+\frac{M'_\nu(z)}{M_\nu(z)}-\theta_{\nu}'(z)\tan\theta_\nu(z).
\]
Second, taking the logarithmic derivative this time in~\eqref{eq:M_id} results in the equality
\[
	\frac{M'_\nu(z)}{M_\nu(z)}=\frac{1}{2z}+\frac{\jj{\nu}{z}\djj{\nu}{z}+\yy{\nu}{z}\dyy{\nu}{z}}{J^2_\nu(z)+Y^2_\nu(z)}.
\]
Combining these two equations, we arrive at~\eqref{eq:log_der_J} with
\[
 \mathcal{Z}_{0}(z;\nu):=\frac{\jj{\nu}{z}\djj{\nu}{z}+\yy{\nu}{z}\dyy{\nu}{z}}{J^2_\nu(z)+Y^2_\nu(z)}.
\]
Next, we apply Lemma~\ref{lem:HankelAsymp} with $p=0$ to the nominator term, with $p=1$ to the denominator term, and estimate as follows,
\[
|\mathcal{Z}_{0}(z;\nu)|=\left|\frac{\mathcal{X}_{0}(z;\nu)}{\frac{2}{\pi z}+\mathcal{Y}_{1}(z;\nu)}\right|\leq\frac{C_{0}\pi}{|z|}\frac{1}{2-\pi C_{1}'\frac{\nu^{2}}{|z|^{2}}}<\frac{C_{0}\pi}{|z|},
\]
where we used that $\pi C_{1}'\nu^{2}/|z|^{2}<1$ for $z\in\mathcal{M}_{\nu}$, see definition~\eqref{eq:set_M}.

Further, we prove~\eqref{eq:log_der_H}.
By using the connection formula~\cite[Eq.~(10.4.3)]{dlmf}
\[
 \hh{\nu}{z}=J_{\nu}(z)+\ii Y_\nu(z),
\]
we get
\[
\frac{\dhh{\nu}{z}}{\hh{\nu}{z}}=\dfrac{\djj{\nu}{z}+\ii\dyy{\nu}{z}}{\jj{\nu}{z}+\ii\yy{\nu}{z}}
	=\frac{\jj{\nu}{z}\djj{\nu}{z}+\yy{\nu}{z}\dyy{\nu}{z}+\ii(\jj\nu{z}\dyy\nu{z}-\djj\nu{z}\yy\nu{z})}{J^2_\nu(z)+Y^2_\nu(z)}.
\]
Employing also the Wronskian identity~\eqref{eq:wronsk_JY} and taking the reciprocal, we arrive at the equality
\begin{equation}
\frac{\hh{\nu}{z}}{\dhh{\nu}{z}}=\frac{J^2_\nu(z)+Y^2_\nu(z)}{\jj{\nu}{z}\djj{\nu}{z}+\yy{\nu}{z}\dyy{\nu}{z}+2\ii/(\pi z)}.
\label{eq:hankel_log_der_inrpoof}
\end{equation}
Taking~\eqref{eq:phase-der_id} into account and applying formula~\eqref{eq:HankelAsymp1} with $p=1$, we get
\[
\theta_{\nu}'(z)\frac{\hh{\nu}{z}}{\dhh{\nu}{z}}=\frac{1}{\ii+\frac{\pi z}{2}\left(\jj\nu{z}\djj\nu{z}+\yy\nu{z}\dyy\nu{z}\right)}=\frac{1}{\ii-\frac{1}{2z}+\frac{\pi z}{2}\mathcal{X}_{1}(z;\nu)}.
\]
From here, one infers~\eqref{eq:log_der_H} with
\[
 \mathcal{Z}_{1}(z;\nu):=\frac{-\frac{1}{z}+\pi z\mathcal{X}_{1}(z;\nu)}{2+\frac{\ii}{z}-\ii\pi z\mathcal{X}_{1}(z;\nu)}.
\]
Estimates similar to those from the first part of the proof yield
\[
|\mathcal{Z}_{1}(z;\nu)|\leq\frac{1}{|z|}\frac{1+\frac{C_{1}}{C_{1}'}}{2-\frac{1}{|z|}\left(1+\frac{C_{1}}{C_{1}'}\right)}<\left(1+\frac{C_1}{C_1'}\right)\frac{1}{|z|},
\]
where we used the additional assumption $|z|>1+C_{1}/C_{1}'$. The proof is complete.
\end{proof}

\begin{remark}\label{rem:bounded_ratio_H}
 Recall that, for $\nu>-1$, the Bessel function $J_{\nu}$ has infinitely many positive simple zeros, see, for example~\cite[\S 10.21(i)]{lozier2003nist}. At these points, it is the first term on the right-hand side of~\eqref{eq:log_der_J} that explodes. The second term represented by $\mathcal{Z}_{0}$ remains bounded. Further, similarly as~\eqref{eq:log_der_H} was proved, one can apply Lemma~\ref{lem:HankelAsymp} in~\eqref{eq:hankel_log_der_inrpoof} to verify that
the ratio $\hh{\nu}{z}/\dhh{\nu}{z}$ remains bounded for all $z\in\mathcal{M}_{\nu}$ with $|z|$ sufficiently large. It follows that the derivative of the analytic function $H_{\nu}^{(1)}$ has no zeros in this region since $\hh{\nu}{z}$ and $\dhh{\nu}{z}$ have no zeros in common which is a consequence of the Wronskian identity 
\[
 J_{\nu}(z)(H_{\nu}^{(1)})'(z)-J_{\nu}'(z)H_{\nu}^{(1)}(z)=\frac{2\ii}{\pi z},
\]
see~\cite[Eqs.~(10.5.3) and~(10.6.2)]{dlmf}.
For more information about the location of zeros of $H_{\nu}^{(1)}$ and its derivative, see~\cite{cruz1982zeros, ker-sko_85}.
\end{remark}

\section{The Demuth--Hansmann--Katriel open problem}\label{sec:OnTheProb}

In this section, we prove Theorem~\ref{ThmMainthm}.

\subsection{Initial set-up}

Rather than analysing the asymptotic behaviour of the entire discrete spectrum of $H_{h}$ for $h$ large, we impose additional restrictions on the discrete eigenvalues and analyse only those from a certain complex set evolving with $h\to\infty$. This is inspired by the approach used in the one-dimensional situation~\cite{bogli2021lieb}. The eigenvalues will be indexed by two indices from certain integer sets for which we need to introduce a suitable notation. 

We fix three parameters 
\begin{equation}
0<\alpha<\beta<\gamma<\frac{1}{2},
\label{eq:param}
\end{equation}
and define
\begin{equation}\label{eq:restr_ell}
\mathcal{L}(h):=\left\{\ell\in\N \mid h^{\alpha+1/2}\leq\ell\leq h^{\beta+1/2}\right\}
\end{equation}
for $h>0$. Next, we fix $\varepsilon\in(0,1)$, put
\begin{equation}
 q:=\frac{1-\varepsilon}{d-1},
\label{eq:def_q}
\end{equation}
and define the set 
\begin{equation}\label{eq:set_J}
\mathcal{J}(h,\ell):=\left\{j\in\N \;\bigg|\; \ell\log^{q}\ell\leq j\leq h^{\gamma+1/2} \right\}
\end{equation}
for $\ell\in\mathbb{N}$ and $h>0$.

Next we slightly reformulate the inequality~\eqref{eq:main-ineq}. Recall definition~\eqref{eq:defH} of $V_{h}$. We have
\[
	\|V_h\|^p_{L^p}=\int_{\mathbb{R}^d}|V_h(x)|^p\,\dd x=\mu_{d}h^{p},
\]
where $\mu_{d}:=\pi^{d/2}/\Gamma(1+d/2)$ is the volume of the unit ball $B_1(0)$ in $\R^{d}$.
Further, Lemma~\ref{LemmaLocationEigenvalues} tells us that $\sigma_{\dd}(H_h)\subset[0,\infty)+\ii[0,h]$. 
Consequently, $\text{dist}(\lambda,[0,\infty))=\im\lambda$ for all $\lambda\in\sigma_{\dd}(H_h)$, and hence
\[
	\frac{1}{\|V_h\|^p_{L^p}}\sum_{\lambda\in\sigma_{\dd}(H_h)}\dfrac{\text{dist}(\lambda,[0,\infty))^p}{|\lambda|^{d/2}}=\frac{1}{\mu_{d} h^p}\sum_{\lambda\in\sigma_{\dd}(H_h)}\dfrac{(\im\lambda)^p}{|\lambda|^{d/2}}.
\]
Therefore, to establish Theorem \ref{ThmMainthm}, it is sufficient to prove the following statement
and then to set $C_{p,d}:=C_{p,d}'\,(\beta-\alpha)$.

\begin{Proposition}\label{LemmaafterMain}
Suppose $d\geq 2$ and $p>0$. Let $\alpha,\beta$ satisfy \eqref{eq:param} and let $0<\varepsilon<1$.
Then there exist $C_{p,d}'>0$ and $h_0>0$ such that, for all $h\geq h_0$, we have
$$\frac{1}{h^p}\sum_{\lambda\in\sigma_{\dd}(H_h)}\dfrac{(\im\lambda)^p}{|\lambda|^{d/2}}\geq 
C_{p,d}'\,(\beta-\alpha)\,(\log h)^{\varepsilon}.$$
\end{Proposition}

\begin{remark}
While the constant $C_{p,d}'$ does not depend on the choice of $\varepsilon$ and $\alpha,\beta,\gamma$, the threshold value $h_0$ depends on these parameters. In fact, keeping track of the error terms below, one sees that they depend on $h^{-\alpha}$, $h^{\gamma-1/2}$ and $(\log h)^{-q}$, which all go to zero as $h\to\infty$ but the convergence rate is dependent on the choice of parameters.
\end{remark}

The proof of Proposition~\ref{LemmaafterMain} proceeds in  five steps worked out in the subsections below.

Further we make the following notational conventions. It would become cumbersome to keep track on all constants in various estimations to be made (as we kept doing in Section~\ref{sec:Bessel}). Therefore we will occasionally use the notation $\lesssim$ or $\gtrsim$ when the inequalities $\leq$ or $\geq$ hold up to a multiplicative constant. Moreover, it is important that the hidden constant will be independent on any choice of indices $\ell\in\mathcal{L}(h)$, $j\in\mathcal{J}(h,\ell)$ (also on the parameters $\alpha,\beta,\gamma$ and $\varepsilon$), and sometimes also on a~complex variable $m$ from a~specified subset of $\C$ depending on $j$ and $\ell$. We will stress this uniformness explicitly whenever convenient. 
For instance, we may say that there exists $h_0>0$ such that for all $h\geq h_{0}$ and $\ell\in\mathcal{L}(h)$, $\nu$ defined by~\eqref{eq:Nudefinition} satisfies $\nu\lesssim\ell$.

In a similar spirit, we will occasionally use the Landau symbol $\mathcal{O}$ for $h\to\infty$, where the hidden constant is independent of $\ell\in\mathcal{L}(h)$, $j\in\mathcal{J}(h,\ell)$ (also of $\alpha,\beta,\gamma$ and $\varepsilon$), and sometimes also of a complex variable $m$ from a  specified subset of $\C$. Nevertheless, we will repeat this uniformness several times and sometimes express the inequality more explicitly for clarity, see for example~\eqref{eq:unif_lim_inproof} or \eqref{eq:err_small}. Lastly, we sometimes, for brevity, say that, for example, a function $f=f(m)$ tends to, say, $1$ as $h\to\infty$ uniformly in $m\in X(\ell,j)\subset\C$, $j\in\mathcal{J}(h,\ell)$, and $\ell\in\mathcal{L}(h)$, by which we mean that
\[
 \lim_{h\to\infty}\sup\{|f(m)-1| \mid m\in X(\ell,j)\subset\C, j\in\mathcal{J}(h,\ell), \ell\in\mathcal{L}(h)\}=0,
\]
or in other words,
\[
 \lim_{h\to\infty}\;\sup_{\ell\in\mathcal{L}(h)}\;\sup_{j\in\mathcal{J}(h,\ell)}\;\sup_{m\in X(\ell,j)}\;|f(m)-1|=0.
\]
Again, in any case of a possible confusion, we express the limit relation in its full.

\subsection{Step~1: Auxiliary zeros}

For $j\in\N$, $\nu>0$, and $h>0$, we define the auxiliary function
\begin{equation}
 f_{\nu,j}(m):=\theta_{\nu}(m)-\frac{\pi}{4}-2\pi j- \ii\log\frac{\sqrt{h}}{4\pi j}
\label{eq:def_f_nu,j}
\end{equation}
and study its zeros in open balls $B_{\nu}(m_{\nu,j}^{(0)})$ of radius $\nu$ centred at points
\begin{equation}
 m_{\nu,j}^{(0)}:=2\pi j+\frac{\nu\pi}{2}+\frac{\pi}{2}+\ii\log\frac{\sqrt{h}}{4\pi j}.
\label{eq:def_m_j^0}
\end{equation}
We designate the dependence on $\nu$ for a clear reference to the order although $\nu$ will be immediately supposed to be determined by equation~\eqref{eq:Nudefinition}, and hence, rather than an independent variable, it is determined by the index $\ell$.

\begin{lemma}\label{lem:aux_zeros}
Let $\nu=\ell+\frac{d}{2}-1$. Then there exists $h_0>0$ such that for all $h\geq h_0$, all $\ell\in\mathcal{L}(h)$ and all $j\in\mathcal{J}(h,\ell)$, the following claims hold true:
\begin{enumerate}[(i)]
\item The function $f_{\nu,j}$ is analytic in the ball $B_{\nu}(m_{\nu,j}^{(0)})$ with a unique simple zero $m_{\nu,j}^{(1)}$ therein;
\item $|m_{\nu,j}^{(1)}-m_{\nu,j}^{(0)}|<\nu/2$;
\end{enumerate}
and, in addition, for any two indices $j_{1},j_{2}\in\mathcal{J}(h,\ell)$, $j_1\neq j_2$, we have 
\begin{enumerate}
\item[(iii)] $|m_{\nu,j_1}^{(1)}-m_{\nu,j_2}^{(1)}|>4$.
\end{enumerate}
\end{lemma}

\begin{proof}
Before proving each claim, we observe that there exists a constant $C>0$ such that 
\begin{equation}
 \sup\left\{\frac{\nu}{\Re m} \;\bigg|\; |m-m_{\nu,j}^{(0)}|\leq\nu, j\in\mathcal{J}(h,\ell), \ell\in\mathcal{L}(h) \right\}\leq\frac{C}{\log^{q}h}
\label{eq:unif_lim_inproof}
\end{equation}
for all $h$ sufficiently large. To see that, we first use the definitions~\eqref{eq:Nudefinition}, \eqref{eq:set_J}, and~\eqref{eq:restr_ell} to estimate
\[
\frac{\Re m_{\nu,j}^{(0)}}{\nu}=\frac{2\pi j}{\nu}+\frac{\pi}{2}+\frac{\pi}{2\nu}\geq\frac{j}{\nu}\gtrsim\frac{j}{\ell}\geq\log^{q}\ell\gtrsim\log^{q}h
\]
for all $h$ sufficiently large, where the non-displayed constants are independent of the choices of $\ell\in\mathcal{L}(h)$ and $j\in\mathcal{J}(h,\ell)$. Hence, for any $m$ in the closure of the ball $B_{\nu}(m_{\nu,j}^{(0)})$, we get
\[
\frac{\Re m}{\nu}\geq \frac{\Re m_{\nu,j}^{(0)}}{\nu}-1\gtrsim\log^{q}h.
\]
It follows that $\Re m>0$ and implies \eqref{eq:unif_lim_inproof}.

\emph{Proof of claim~(i):}
We prove the analyticity of $f_{\nu,j}$ in $B_{\nu}(m_{\nu,j}^{(0)})$ by showing that $B_{\nu}(m_{\nu,j}^{(0)})\subset\mathcal{M}_{\nu}$ for all $h$ sufficienly large. Since $\theta_{\nu}$ is analytic in $\mathcal{M}_{\nu}$, see Lemma~\ref{lem:phase_analytic}, the analyticity of $f_{\nu,j}$ in $B_{\nu}(m_{\nu,j}^{(0)})$ then follows immediately from its definition~\eqref{eq:def_f_nu,j}.
 
First, it follows readily from~\eqref{eq:unif_lim_inproof} that, for all $h$ sufficiently large, we have $A\nu<\Re m$ for any fixed $A>0$ (hence in particular for the constant $A$ in~\eqref{eq:defA} used in the definition~\eqref{eq:set_M} of the set $\mathcal{M}_{\nu}$) and all $m\in B_{\nu}(m_{\nu,j}^{(0)})$. Second, recalling definitions~\eqref{eq:def_m_j^0} and \eqref{eq:set_J}, we find
\[
 |\Im m_{\nu,j}^{(0)}|=\log\frac{4\pi j}{\sqrt{h}}\lesssim\log h,
\]
and taking also~\eqref{eq:Nudefinition} and~\eqref{eq:restr_ell} into account, we get
\[
 \frac{|\Im m_{\nu,j}^{(0)}|}{\nu}\lesssim\frac{\log h}{h^{\alpha+1/2}}.
\]
As a result, for $m\in B_{\nu}(m_{\nu,j}^{(0)})$, we deduce that 
\[
 \frac{|\Im m|}{\Re m}\leq\frac{1+|\Im(m_{\nu,j}^{(0)})|/\nu}{-1+\Re(m_{\nu,j}^{(0)})/\nu}\lesssim\frac{1}{\log^{q} h},
\]
which implies that $|m|<2\Re m$ for all $h$ sufficiently large. Recalling the definition~\eqref{eq:set_M} of $\mathcal{M}_{\nu}$, we have shown that, for all $h$ large enough, we have $B_{\nu}(m_{\nu,j}^{(0)})\subset\mathcal{M}_{\nu}$ for all $j\in\mathcal{J}(h,\ell)$ and $\ell\in\mathcal{L}(h)$.
 
Next, we prove that $f_{\nu,j}$ has a unique simple zero in $B_{\nu}(m_{\nu,j}^{(0)})$. Using~\eqref{eq:def_psi}, we may write
\[
 f_{\nu,j}(m)=g_{\nu,j}(m)-\psi_{\nu}(m),
\]
where $g_{\nu,j}(m):=m-m_{j,\nu}^{(0)}$. Clearly, $m_{j,\nu}^{(0)}$ is a unique simple zero of $g_{\nu,j}(m)$ in $\C$, so also in $B_{\nu}(m_{\nu,j}^{(0)})$. Bearing in mind~\eqref{eq:unif_lim_inproof}, from which it follows that $\nu/|m|\to0$ for $h\to\infty$ uniformly in the closure of the ball $B_{\nu}(m_{\nu,j}^{(0)})$, we infer from Lemma~\ref{lem:PhaseAsymp} that 
\[
 |f_{\nu,j}(m)-g_{j,\nu}(m)|=|\psi_{\nu}(m)|<\nu=|g_{j,\nu}(m)|
\]
for all $m\in\partial B_{\nu}(m_{\nu,j}^{(0)})$ and $h$ sufficiently large. By Rouch{\' e}'s theorem, $f_{\nu,j}$ and $g_{\nu,j}$ have the same number of zeros including multiplicities in $B_{\nu}(m_{\nu,j}^{(0)})$, i.e. $f_{\nu,j}$ has a unique simple zero in $B_{\nu}(m_{\nu,j}^{(0)})$, which we denote by $m_{\nu,j}^{(1)}$.

\emph{Proof of claim~(ii):}
From the definition of $f_{\nu,j}$ and Lemma~\ref{lem:PhaseAsymp}, we infer that
\[
|m_{\nu,j}^{(1)}-m_{\nu,j}^{(0)}|=|\psi_{\nu}(m_{\nu,j}^{(1)})|\leq 2A^{2}\frac{\nu^{2}}{|m_{\nu,j}^{(1)}|}<\frac{\nu}{2}
\]
for all $h$ sufficiently large since $\nu/m_{\nu,j}^{(1)}\to0$ as $h\to\infty$ by \eqref{eq:unif_lim_inproof}.

\emph{Proof of claim~(iii):}
Let $j_{1},j_{2}\in\mathcal{J}(h,\ell)$, $j_1\neq j_2$. Then from the definition of $f_{\nu,j}$ and~\eqref{eq:def_psi}, one gets
\[
 m_{\nu,j_1}^{(1)}-m_{\nu,j_2}^{(1)}=2\pi(j_1-j_2)+\psi_{\nu}(m_{\nu,j_1}^{(1)})-\psi_{\nu}(m_{\nu,j_2}^{(1)})+\ii\log\frac{j_2}{j_1}.
\]
It follows the estimate
\begin{equation}
 |m_{\nu,j_1}^{(1)}-m_{\nu,j_2}^{(1)}|\geq 2\pi|j_1-j_2|-|\psi_{\nu}(m_{\nu,j_1}^{(1)})-\psi_{\nu}(m_{\nu,j_2}^{(1)})|.
 \label{eq:diff_m_inproof}
\end{equation}
By the mean value theorem, 
\[
|\psi_{\nu}(m_{\nu,j_1}^{(1)})-\psi_{\nu}(m_{\nu,j_2}^{(1)})|
\leq 
|m_{\nu,j_1}^{(1)} - m_{\nu,j_2}^{(1)}|
\sup_{m\in(m_{\nu,j_1}^{(1)},m_{\nu,j_2}^{(1)})}|\psi_{\nu}'(m)|
\]
where $(m_{\nu,j_1}^{(1)},m_{\nu,j_2}^{(1)})$ is the open complex line segment connecting the points $m_{\nu,j_1}^{(1)}$ and $m_{\nu,j_2}^{(1)}$, which is entirely located in $\mathcal{M}_{\nu}$. Since any $m\in(m_{\nu,j_1}^{(1)},m_{\nu,j_2}^{(1)})$ satisfies
\[
 |m|\geq\min(\Re m_{\nu,j_1}^{(1)},\Re m_{\nu,j_2}^{(1)}),
\]
we deduce, with the aid of~\eqref{eq:psi_der_bound}, the upper bound
\[
\sup_{m\in (m_{\nu,j_1}^{(1)},m_{\nu,j_2}^{(1)})}|\psi_{\nu}'(m)|\leq\frac{A^{2}\nu^{2}}{\min\!\left((\Re m_{\nu,j_1}^{(1)})^{2},(\Re m_{\nu,j_2}^{(1)})^{2}\right)}.
\]
Employing~\eqref{eq:unif_lim_inproof} once more, we can take $h$ sufficiently large to ensure that 
\[
\sup_{m\in (m_{\nu,j_1}^{(1)},m_{\nu,j_2}^{(1)})}|\psi_{\nu}'(m)|<\frac{1}{2},
\]
and hence 
\[
|\psi_{\nu}(m_{\nu,j_1}^{(1)})-\psi_{\nu}(m_{\nu,j_2}^{(1)})|
\leq\frac{1}{2}\,|m_{\nu,j_1}^{(1)} - m_{\nu,j_2}^{(1)}|.
\]
Plugging the last estimate into~\eqref{eq:diff_m_inproof}, we obtain 
\[
 |m_{\nu,j_1}^{(1)}-m_{\nu,j_2}^{(1)}|\geq 2\pi|j_1-j_2|-\frac{1}{2}\,|m_{\nu,j_1}^{(1)} - m_{\nu,j_2}^{(1)}|,
\]
from which, bearing in mind that $j_1\neq j_2$, we conclude
\[
|m_{\nu,j_1}^{(1)}-m_{\nu,j_2}^{(1)}|\geq \frac{4\pi}{3}|j_1-j_2|\geq \frac{4\pi}{3}>4.
\]
The proof is complete.
\end{proof}

Further, we restrict our analysis to fixed neighborhoods of the zeros $m_{\nu,j}^{(1)}$ from Lemma~\ref{lem:aux_zeros}. Concretely, we consider balls $B_{2}(m_{\nu,j}^{(1)})$ of radius $2$ centred at $m_{\nu,j}^{(1)}$. We may always suppose that $h_0$ is large enough so that for all $h\geq h_{0}$ we have
\begin{equation}\label{eq:B2_in_M}
B_{2}(m_{\nu,j}^{(1)})\subset B_{\nu}(m_{\nu,j}^{(0)})\subset\mathcal{M}_{\nu}
\end{equation}
and
\begin{equation}\label{eq:B2_disjoint}
B_{2}(m_{\nu,j_1}^{(1)})\cap B_{2}(m_{\nu,j_2}^{(1)})=\emptyset
\end{equation}
for any $j,j_1,j_2\in\mathcal{J}(h,\ell)$, $j_1\neq j_2$, and $\ell\in\mathcal{L}(h)$. The first inclusion in~\eqref{eq:B2_in_M} is a consequence of claim (ii) of Lemma~\ref{lem:aux_zeros}.
Note that we have used that $\nu\geq 4$ which is satisfied for all sufficiently large $h$ since $\nu=\ell+\frac d 2 -1$ and $\ell\in\mathcal L(h)$ diverges as $h\to\infty$.
The second inclusion in~\eqref{eq:B2_in_M} has been verified in the proof of Lemma~\ref{lem:aux_zeros}. The disjointness~\eqref{eq:B2_disjoint} follows immediately from claim (iii) of Lemma~\ref{lem:aux_zeros}.

\begin{lemma}\label{lem:sup_re-im}
Let $\nu=\ell+\frac{d}{2}-1$. For $h\to\infty$, we have the limits
\begin{equation}
\sup\left\{\left|\frac{\Re m}{2\pi j}-1\right|\;\bigg|\; |m-m_{\nu,j}^{(1)}|\leq2, j\in\mathcal{J}(h,\ell), \ell\in\mathcal{L}(h) \right\}\to0
\label{eq:sup_re_lim}
\end{equation}
and
\begin{equation}
 \sup\left\{\left|\frac{\Im m}{\log\left(\sqrt{h}/(4\pi j)\right)}-1\right|\;\bigg|\; |m-m_{\nu,j}^{(1)}|\leq2, j\in\mathcal{J}(h,\ell), \ell\in\mathcal{L}(h) \right\}\to0.
\label{eq:sup_im_lim}
\end{equation}
\end{lemma}

\begin{proof}
 Notice that, since for $h\to\infty$ we have 
 \[
  1/j\to0 \quad\mbox{ and }\quad 1\big/\log\frac{\sqrt{h}}{4\pi j}\to0
 \]
 uniformly in $j\in\mathcal{J}(h,\ell)$ and $\ell\in\mathcal{L}(h)$, it is sufficient, in order to prove \eqref{eq:sup_re_lim} and \eqref{eq:sup_im_lim}, to verify the limits
\begin{equation}
\lim_{h\to\infty}\sup\left\{\left|\frac{\Re m_{\nu,j}^{(1)}}{2\pi j}-1\right|\;\bigg|\;  j\in\mathcal{J}(h,\ell), \ell\in\mathcal{L}(h)\right\}=0
\label{eq:sup_re_lim_inproof}
\end{equation}
and
\begin{equation}
 \lim_{h\to\infty}\sup\left\{\left|\frac{\Im m_{\nu,j}^{(1)}}{\log\left(\sqrt{h}/(4\pi j)\right)}-1\right|\;\bigg|\; j\in\mathcal{J}(h,\ell), \ell\in\mathcal{L}(h)\right\}=0.
\label{eq:sup_im_lim_inproof}
\end{equation}

First we show~\eqref{eq:sup_re_lim_inproof}. Taking real parts in the equation $f_{\nu,j}(m_{\nu,j}^{(1)})=0$ and recalling~\eqref{eq:def_psi}, 
one gets
\begin{equation}
 \Re m_{\nu,j}^{(1)}=2\pi j+\frac{\nu\pi}{2}+\frac{\pi}{2}+\Re \psi_{\nu}(m_{\nu,j}^{(1)}).
\label{eq:re_m_1}
\end{equation}
For $h\to\infty$, the asymptotically dominating term on the right hand side is $2\pi j$. Indeed, we infer from definitions~\eqref{eq:Nudefinition}, \eqref{eq:set_J}, and \eqref{eq:restr_ell} that $\nu/j\lesssim\log^{-q}h$. Taking also into account that $|\psi_{\nu}(m_{\nu,j}^{(1)})|<\nu$, as it follows from~\eqref{eq:PhaseAsymp} for all $h$ large, and bearing in mind that $m_{\nu,j}^{(1)}\in\mathcal{M}_{\nu}$, see~\eqref{eq:B2_in_M}, we verify~\eqref{eq:sup_re_lim_inproof}.

In order to check the limit~\eqref{eq:sup_im_lim_inproof}, we proceed similarly by taking imaginary parts in $f_{\nu,j}(m_{\nu,j}^{(1)})=0$. It  implies the identity 
\begin{equation}
\Im m_{\nu,j}^{(1)}=\Im \psi_{\nu}(m_{\nu,j}^{(1)})+\log\frac{\sqrt{h}}{4\pi j}.
\label{eq:im_m_1}
\end{equation}
One infers from~\eqref{eq:imaginaryPhase} that, uniformly in $\ell\in\mathcal{L}(h)$ and $j\in\mathcal{J}(h,\ell)$, we have
\[
\frac{\Im \psi_{\nu}(m_{\nu,j}^{(1)})}{\Im m_{\nu,j}^{(1)}}\to0
\]
since $\nu/m_{\nu,j}^{(1)}\to0$ as $h\to\infty$, see~\eqref{eq:unif_lim_inproof}. The formula~\eqref{eq:sup_im_lim_inproof} follows.
\end{proof}

\begin{remark}
 Notice that the proof of~\eqref{eq:sup_re_lim} actually shows more on the decay of the remainder. Namely, from equality~\eqref{eq:re_m_1} and estimates made in the paragraph below, it follows that there is $h_0>0$ such that for all $h\geq h_{0}$, $\ell\in\mathcal{L}(h)$, $j\in\mathcal{J}(h,\ell)$, and  $m\in B_{2}(m_{\nu,j}^{(1)})$, we have the estimate 
\begin{equation}\label{eq:re_m_reminder}
\left|\frac{\Re m}{2\pi j}-1\right|\leq\frac{C}{\log^{q}h},
\end{equation}
where $C>0$ is a constant that is independent of $\ell$, $j$, and $m$.
\end{remark}

\begin{remark}\label{rem:fourth-quadrant}
 It follows from formulas~\eqref{eq:sup_re_lim} and~\eqref{eq:sup_im_lim} that, for all $h>0$ sufficiently large, the balls $B_{2}(m_{\nu,j}^{(1)})$ are located in the fourth quadrant of the complex plane ($\Re m>0$ and $\Im m<0$) for all $\ell\in\mathcal{L}(h)$ and $j\in\mathcal{J}(h,\ell)$.
\end{remark}

Next, we pass to locating the parameter $k=k(m)$ determined by $m\in B_{2}(m_{\nu,j}^{(1)})$ and the formula
\begin{equation}\label{eq:defk}
k=\sqrt{\ii h +m^{2}},
\end{equation}
where the square root assumes its principal branch, see~\eqref{eq:settingLambda}.

\begin{lemma}\label{lem:param_k}
Let $\nu=\ell+\frac{d}{2}-1$. There is $h_{0}>0$ such that, for all $h\geq h_{0}$, if $m\in B_{2}(m_{\nu,j}^{(1)})$ for any $\ell\in\mathcal{L}(h)$ and $j\in\mathcal{J}(h,\ell)$, then $k\in\mathcal{M}_{\nu}$. Moreover, as $h\to\infty$, we have the limits 
\begin{equation}
 \frac k m \to 1
 \quad\mbox{ and }\quad
\frac{k}{\Re m}\to1
\label{eq:lim_abs-k}
\end{equation}
uniformly in $m\in B_{2}(m_{\nu,j}^{(1)})$, $j\in\mathcal{J}(h,\ell)$, and $\ell\in\mathcal{L}(h)$.
\end{lemma}

\begin{proof}
First, we deduce the uniform limits from~\eqref{eq:lim_abs-k} and prepare further auxiliary limits for the claim $k\in\mathcal{M}_{\nu}$, which is to be proven afterwards.
It follows from the equality $k^{2}=\ii h+m^{2}$ that
\begin{equation}\label{eq:k^2}
 \begin{aligned}
	k^2&=(\rr m)^2-(\im m)^2+\ii\left(h+2\rr(m)\im(m)\right)\\
	&=(\rr m)^2\left[1-\frac{(\im m)^2}{(\rr m)^2}+\ii\left(\frac{h}{(\rr m)^2}+\frac{2\im m }{\rr m}\right)\right].
\end{aligned}
\end{equation}
As $h\to\infty$, we have 
\begin{equation}\label{eq:zerolimits}
\frac{\im m}{\rr m}\to0
\quad\mbox{ and }\quad
\frac{h}{(\rr m)^2}\to0
\end{equation}
uniformly in $m\in B_{2}(m_{\nu,j}^{(1)})$, $j\in\mathcal{J}(h,\ell)$, and $\ell\in\mathcal{L}(h)$. The first zero-limit formula follows by inspection of the limits from Lemma~\ref{lem:sup_re-im}. To verify the second zero-limit, we estimate
\[
 \frac{h}{(\rr m)^2}\lesssim\frac{h}{j^{2}}\leq\frac{h}{\ell^2}\leq\frac{1}{h^{2\alpha}}
\]
using the ranges for indices $\ell$ and $j$ from definitions~\eqref{eq:restr_ell} and~\eqref{eq:set_J}. It follows that
\begin{equation}
 \frac{\Re k^{2}}{(\Re m)^{2}}\to 1 
 \quad\mbox{ and }\quad
 \frac{\Im k^{2}}{(\Re m)^{2}}\to 0,
 \label{eq:lim_re-k2_inproof}
\end{equation}
and so
\begin{equation}
\frac{|k|^{2}}{(\Re m)^{2}}\to 1,
\label{eq:lim_abs-k2_inproof}
\end{equation}
as $h\to\infty$, uniformly in $m\in B_{2}(m_{\nu,j}^{(1)})$, $j\in\mathcal{J}(h,\ell)$, and $\ell\in\mathcal{L}(h)$. 
Since $k$ is the principal branch square root of the expression in \eqref{eq:k^2}, using \eqref{eq:zerolimits}
and $\Re m>0$ by Remark~\ref{rem:fourth-quadrant}, we obtain the second limit in~\eqref{eq:lim_abs-k}. 
Further, we infer from the first limit in~\eqref{eq:zerolimits} that 
\[
 \frac{\Re m}{m}\to1.
\]
Consequently, the first limit in~\eqref{eq:lim_abs-k} follows from the second one.

By~\eqref{eq:set_M}, the parameter $k=k(m)$ is in $\mathcal{M}_{\nu}$, if $A\nu<\Re k$ and $|k|<2\Re k$. We show that this is the case for all $h$ sufficiently large independently on a choice of $m\in B_{2}(m_{\nu,j}^{(1)})$, $j\in\mathcal{J}(h,\ell)$, and $\ell\in\mathcal{L}(h)$. First, we verify the inequality $|k|<2\Re k$. Since $k$ is defined as the pricipal branch square root of a number, which is not purely negative, we know that $\Re k>0$. Therefore the inequality $|k|<2\Re k$ holds true if
\begin{equation}
 \frac{2(\Re k)^{2}}{|k|^{2}}>\frac{1}{2}.
\label{eq:k_in_M_first_inproof}
\end{equation}
Using the identity 
\[
 \Re k^{2}=(\Re k)^{2}-(\Im k)^{2}=2(\Re k)^{2}-|k|^{2},
\]
we find that 
\[
\frac{2(\Re k)^{2}}{|k|^{2}}=1+\frac{\Re k^{2}}{|k|^{2}}=1+\frac{\Re k^{2}}{(\Re m)^{2}}\frac{(\Re m)^{2}}{|k|^{2}}\to 2
\]
as $h\to\infty$, uniformly in $m\in B_{2}(m_{\nu,j}^{(1)})$, $j\in\mathcal{J}(h,\ell)$, and $\ell\in\mathcal{L}(h)$, which follows from the limits in~\eqref{eq:lim_re-k2_inproof} and~\eqref{eq:lim_abs-k2_inproof}.
This yields~\eqref{eq:k_in_M_first_inproof} for all $h$ large enough.

Next, we write
\[
 \frac{\nu^{2}}{(\Re k)^{2}}=\frac{\nu^{2}}{(\Re m)^{2}}\frac{(\Re m)^{2}}{|k|^{2}}\frac{|k|^{2}}{(\Re k)^{2}}.
\]
We already know that that the fractions
\[
\frac{(\Re m)^{2}}{|k|^{2}}
\quad\mbox{ and }\quad
\frac{|k|^{2}}{(\Re k)^{2}}
\]
converge to $1$ as $h\to\infty$. Moreover, recalling \eqref{eq:unif_lim_inproof}, we find that
\[
\frac{\nu^{2}}{(\Re m)^{2}}\to0
\]
as $h\to\infty$, uniformly in $m\in B_{2}(m_{\nu,j}^{(1)})$, $j\in\mathcal{J}(h,\ell)$, and $\ell\in\mathcal{L}(h)$. In total, we observe that 
\[
 \frac{\nu^{2}}{(\Re k)^{2}}\to0
\]
uniformly as $h\to\infty$, and this limit, together with the positivity of $\Re k$, yields the inequality $A\nu<\Re k$ for all $h$ sufficiently large. The proof is complete.
\end{proof}

\subsection{Step 2: An auxiliary error function}

First we deduce a few auxiliary inequalities.

\begin{lemma}
Let $\nu=\ell+\frac{d}{2}-1$ and  let $\alpha,\gamma$ satisfy \eqref{eq:param}. Then there exists $h_0>0$ such that for all $h\geq h_0$, all $\ell\in\mathcal{L}(h)$, $j\in\mathcal{J}(h,\ell)$, and $m\in B_{2}(m_{\nu,j}^{(1)})$, we have the inequalities
 \begin{equation}\label{eq:im_theta_bounds}
  -2\gamma\log h\leq\im\theta_{\nu}(m)\leq-\alpha\log h
 \end{equation}
and
 \begin{equation}\label{eq:tan-cos_theta_bounds}
 |\ii+\tan\theta_{\nu}(m)|\leq 4 h^{-2\alpha}
 \quad\mbox{ and }\quad
 |\cos\theta_{\nu}(m)|\leq C\, \frac{j}{\sqrt{h}}
 \end{equation}
 with constant $C:=4\pi \e^4>0$.
\end{lemma}

\begin{proof}
First we establish~\eqref{eq:im_theta_bounds}.
 Recalling the uniform limit~\eqref{eq:unif_lim_inproof}, we may suppose $h_{0}$ to be chosen large enough so that 
 \[
  A^{2}\frac{\nu^{2}}{|m|^{2}}<1
 \]
 for all $m\in B_{2}(m_{\nu,j}^{(1)})$, $j\in\mathcal{J}(h,\ell)$, $\ell\in\mathcal{L}(h)$, and $h\geq h_0$.
 With this choice and by the mean value theorem, we get
 \begin{equation}
  |\theta_{\nu}(m)-\theta_{\nu}(m_{\nu,j}^{(1)})|\leq2\sup_{z\in B_{2}(m_{\nu,j}^{(1)})}|\theta_{\nu}'(z)|\leq 2\bigg(1+\sup_{z\in B_{2}(m_{\nu,j}^{(1)})}|\psi_{\nu}'(z)|\bigg)<4
 \label{eq:theta-increm_inproof}
 \end{equation}
 for all $m\in B_{2}(m_{\nu,j}^{(1)})$, where definition~\eqref{eq:def_psi}, 
and the uniform bound~\eqref{eq:psi_der_bound} were used. Recalling the equality~\eqref{eq:im_m_1}, we have
 \begin{equation}\label{eq:Imtheta}
 \Im\theta_{\nu}(m_{\nu,j}^{(1)})=\log\frac{\sqrt{h}}{4\pi j}.
 \end{equation}
 Together with~\eqref{eq:theta-increm_inproof}, it follows that
 \[
 \left|\im\theta_{\nu}(m)-\log\frac{\sqrt{h}}{4\pi j}\right|<4
 \]
 for any $m\in B_{2}(m_{\nu,j}^{(1)})$. By using the ranges from definitions~\eqref{eq:set_J} and~\eqref{eq:restr_ell} of the sets $\mathcal{J}(h,\ell)$ and $\ell\in\mathcal{L}(h)$, we deduce the two-sided estimate
 \[
 -\log 4\pi-\gamma\log h\leq\log\frac{\sqrt{h}}{4\pi j}\leq-\alpha\log h-\log\log^{q}h-q\log\left(\alpha+\frac{1}{2}\right).
 \]
 Clearly, choosing $h_0$ large enough, we have the inequalities~\eqref{eq:im_theta_bounds}.
 
 To deduce the inequality for the tangent in~\eqref{eq:tan-cos_theta_bounds}, we write
 \[
  \ii+\tan\theta_{\nu}(m)=\ii-\ii\frac{1-\e^{-2\ii\theta_{\nu}(m)}}{1+\e^{-2\ii\theta_{\nu}(m)}}=\frac{2\ii\e^{-2\ii\theta_{\nu}(m)}}{1+\e^{-2\ii\theta_{\nu}(m)}}.
 \]
With the aid of~\eqref{eq:im_theta_bounds}, for $h_0$ sufficiently large, we get the desired estimate
 \[
 |\ii+\tan\theta_{\nu}(m)|\leq\frac{2\e^{2\Im\theta_{\nu}(m)}}{1-\e^{2\Im\theta_{\nu}(m)}}\leq\frac{2h^{-2\alpha}}{1-h^{-2\alpha}}\leq 4 h^{-2\alpha}.
 \]
 
 It remains to verify the inequality for the cosine in~\eqref{eq:tan-cos_theta_bounds}. Using~\eqref{eq:im_theta_bounds}, we find, again for $h_0$ sufficiently large,
 \[
  |\cos\theta_{\nu}(m)|\leq\frac{\e^{-\Im\theta_{\nu}(m)}}{2}\left(1+\e^{2\Im\theta_{\nu}(m)}\right)\leq\frac{\e^{-\Im\theta_{\nu}(m)}}{2}\left(1+h^{-2\alpha}\right)\leq \e^{-\Im\theta_{\nu}(m)}.
 \]
By using also estimate~\eqref{eq:theta-increm_inproof}, we get
\[
   |\cos\theta_{\nu}(m)|\leq \e^{4-\Im\theta_{\nu}(m_{\nu,j}^{(1)})}. 
\]
Finally, with \eqref{eq:Imtheta}, we arrive at the second estimate from~\eqref{eq:tan-cos_theta_bounds}.
\end{proof}

Next we explore properties of an error term function
\begin{equation}\label{eq:def_xi}
 \xi_{\nu}(m):=\sin^{2}\theta_{\nu}(m)+\left(\frac{J_{\nu}'(m)H_{\nu}^{(1)}(k)}{J_{\nu}(m)(H_{\nu}^{(1)})'(k)}\right)^{\!2}\!\cos^{2}\theta_{\nu}(m),
\end{equation}
where $k=\sqrt{\ii h+m^{2}}$ with the principal branch of the square root.

\begin{lemma}\label{lem:xi_prop}
Let $\nu=\ell+\frac{d}{2}-1$ and  $q=\frac{1-\varepsilon}{d-1}$ for a fixed $\varepsilon\in (0,1)$. Then there exists $h_0>0$ such that for all $h\geq h_0$, $\ell\in\mathcal{L}(h)$, and $j\in\mathcal{J}(h,\ell)$, the function $\xi_{\nu}$ is analytic in $B_{2}(m_{\nu,j}^{(1)})$ and there is a~constant $C>0$ independent of $j$, $\ell$, and $m$ such that 
 \begin{equation}\label{eq:xi_small}
  |\xi_{\nu}(m)|\leq\frac{C}{\log^{2q}h}
 \end{equation}
 for any $m\in B_{2}(m_{\nu,j}^{(1)})$.
\end{lemma}

\begin{proof}
First we verify the analyticity of $\xi_{\nu}$ in the balls $B_{2}(m_{\nu,j}^{(1)})$. We may suppose $h$ to be sufficiently large so that $B_{2}(m_{\nu,j}^{(1)})\subset\mathcal{M}_{\nu}$, see~\eqref{eq:B2_in_M}. Then the phase function $\theta_{\nu}$, which appears in definition~\eqref{eq:def_xi} of the function $\xi_{\nu}$, is analytic in~$B_{2}(m_{\nu,j}^{(1)})$ by Lemma~\ref{lem:phase_analytic}. Recall that the Bessel functions $J_{\nu}$ as well as $H_{\nu}^{(1)}$ are analytic functions of their main argument in the right half-plane. The balls $B_{2}(m_{\nu,j}^{(1)})$ are located therein, in fact in the fourth quadrant of $\C$ by Remark~\ref{rem:fourth-quadrant}. The argument $k$, as the principal branch square root of $\ii h+m^{2}$, also fulfills $\Re k>0$. 

Thus, we see from the definition~\eqref{eq:def_xi} that $\xi_{\nu}$ is analytic in $B_{2}(m_{\nu,j}^{(1)})$ if $J_{\nu}(m)$ and $(H_{\nu}^{(1)})'(k)$ do not vanish in $B_{2}(m_{\nu,j}^{(1)})$. The case of $J_{\nu}(m)$ is clear as it is well known that for $\nu>-1$, the Bessel function $J_{\nu}$ possesses real zeros only, see Remark~\ref{rem:bounded_ratio_H}. The same remark also implies that $(H_{\nu}^{(1)})'(k)\neq0$ for all $h$ sufficiently large since $k\in\mathcal{M}_{\nu}$ by Lemma~\ref{lem:param_k}. Consequently, we see that for all $h$ large enough, the function $\xi_{\nu}$ defined by the expression~\eqref{eq:def_xi} is analytic in $B_{2}(m_{\nu,j}^{(1)})$ for all $\ell\in\mathcal{L}(h)$ and $j\in\mathcal{J}(h,\ell)$.

In the rest of the proof, we derive the uniform bound~\eqref{eq:xi_small}. The proof relies on asymptotic analysis of the Bessel functions appearing in the definition~\eqref{eq:def_xi}. Below we make use of the Landau symbol $\mathcal{O}$ for $h\to\infty$ which is uniform in $\ell$, $j$ and $m$, i.e. the involved constant is always independent of $\ell\in\mathcal{L}(h)$, $j\in\mathcal{J}(h,\ell)$, and $m\in B_{2}(m_{\nu,j}^{(1)})$.

Let us rewrite~\eqref{eq:def_xi} as
\begin{equation}
\xi_{\nu}(m)=\left[\tan^{2}\theta_{\nu}(m)+\left(\frac{J_{\nu}'(m)H_{\nu}^{(1)}(k)}{J_{\nu}(m)(H_{\nu}^{(1)})'(k)}\right)^{\!2}\,\right]\cos^{2}\theta_{\nu}(m).
\label{eq:xi_inproof}
\end{equation}
Notice that $\cos\theta_{\nu}(m)\neq0$ for $m\in B_{2}(m_{\nu,j}^{(1)})$ since our $m$ is non-real by Remark~\ref{rem:fourth-quadrant}, while $\cos\theta_{\nu}(m)$ may vanish only on zeros of $J_{\nu}$ that are real, see~\eqref{eq:J-Y_M-psi}.
First we analyse the asymptotic behaviour of the expression from~\eqref{eq:xi_inproof} in the square brackets.
It follows readily from Lemma~\ref{lem:log_der} that
\[
  \frac{\djj{\nu}{m}}{\jj{\nu}{m}}=-\theta_{\nu}'(m)\tan\theta_\nu(m)+\bigO{\frac{1}{m}}
\]
as $h\to\infty$. Using also that the ratio $H_{\nu}^{(1)}(k)/(H_{\nu}^{(1)})'(k)$ remains bounded for $k\in\mathcal{M}_{\nu}$ of sufficiently large modulus, see Remark~\ref{rem:bounded_ratio_H}, we find
\begin{equation}\label{eq:aux_asympt0_inproof}
  \frac{J_{\nu}'(m)H_{\nu}^{(1)}(k)}{J_{\nu}(m)(H_{\nu}^{(1)})'(k)}=-\theta_{\nu}'(m)\frac{H_{\nu}^{(1)}(k)}{(H_{\nu}^{(1)})'(k)}\tan\theta_\nu(m)+\bigO{\frac{1}{m}}.
\end{equation}
The function in front of $\tan\theta_\nu(m)$ cannot be estimated directly by using~\eqref{eq:log_der_H} as the arguments are not the same. Next, we slightly manipulate the expression to deduce its asymptotic behaviour.

By the mean value theorem,
\[
 |\theta_{\nu}'(m)-\theta_{\nu}'(k)|\leq|m-k|\max_{z\in(m,k)}|\theta_{\nu}''(z)|.
\]
For $h$ large enough, the entire line segment $(m,k)$ is located in the convex set $\mathcal{M}_{\nu}$, and so we may apply the uniform bound from~\eqref{eq:psi_der2_bound}. 
Recalling that $\psi_{\nu}''=-\theta_{\nu}''$ and taking~\eqref{eq:lim_abs-k} into account, we obtain
\[
 \max_{z\in(m,k)}|\theta_{\nu}''(z)|=\bigO{\frac{\nu^{2}}{m^{3}}}.
\]
Since $k^{2}-m^{2}=\ii h$ and using~\eqref{eq:lim_abs-k} once more, we find
\[
 |m-k|=\frac{h}{|m+k|}=\bigO{\frac{h}{m}}.
\]
Altogether, we observe that
\[
 \theta_{\nu}'(m)-\theta_{\nu}'(k)=\bigO{\frac{h\nu^{2}}{m^{4}}}.
\]
Now we may estimate 
\[
 \theta_{\nu}'(m)\frac{H_{\nu}^{(1)}(k)}{(H_{\nu}^{(1)})'(k)}=\theta_{\nu}'(k)\frac{H_{\nu}^{(1)}(k)}{(H_{\nu}^{(1)})'(k)}+\left(\theta_{\nu}'(m)-\theta_{\nu}'(k)\right)\frac{H_{\nu}^{(1)}(k)}{(H_{\nu}^{(1)})'(k)}
\]
with the aid of~\eqref{eq:log_der_H} and the boundedness of $H_{\nu}^{(1)}(k)/(H_{\nu}^{(1)})'(k)$ to deduce the uniform asymptotic formula
\begin{equation}\label{eq:aux_asympt1_inproof}
 \theta_{\nu}'(m)\frac{H_{\nu}^{(1)}(k)}{(H_{\nu}^{(1)})'(k)}=-\ii+\bigO{\frac{1}{m}}+\bigO{\frac{h\nu^{2}}{m^{4}}}.
\end{equation}

Plugging~\eqref{eq:aux_asympt1_inproof} into \eqref{eq:aux_asympt0_inproof} and taking also into account the uniform boundedness of $\tan\theta_{\nu}(m)$ in $B_{2}(m_{\nu,j}^{(1)})$, see~\eqref{eq:tan-cos_theta_bounds}, we arrive at the uniform asymptotic formula
\begin{equation}\label{eq:BesselProductAsymp}
  \frac{J_{\nu}'(m)H_{\nu}^{(1)}(k)}{J_{\nu}(m)(H_{\nu}^{(1)})'(k)}=\ii\tan\theta_\nu(m)+\bigO{\frac{1}{m}}+\bigO{\frac{h\nu^{2}}{m^{4}}}.
\end{equation}
When applied in~\eqref{eq:xi_inproof} and using the first formula from~\eqref{eq:tan-cos_theta_bounds} once more, we get
\[
 \xi_{\nu}(m)=\left[\bigO{\frac{1}{m}}+\bigO{\frac{h\nu^{2}}{m^{4}}}\right]\cos^{2}\theta_{\nu}(m).
\]
Using also the second formula from~\eqref{eq:tan-cos_theta_bounds}, we finally deduce that 
\[
 \xi_{\nu}(m)=\bigO{\frac{j^{2}}{mh}}+\bigO{\frac{j^{2}\nu^{2}}{m^{4}}}.
\]

With the aid of Lemma~\ref{lem:sup_re-im} we may estimate the last two uniform Landau $\mathcal{O}s$ as follows. The first error term 
can be estimated by
\[
\frac{j^{2}}{|m|h}\lesssim\frac{j}{h}\leq h^{\gamma-1/2},
\]
whereas the second inequality uses~\eqref{eq:set_J}; note that $\gamma<1/2$. On the other hand, the rate of the second term fulfills
\[
\frac{j^{2}\nu^{2}}{|m|^{4}}\lesssim\frac{\nu^{2}}{|m|^{2}}\lesssim\frac{1}{\log^{2q}h},
\]
where the second inequality is a consequence of~\eqref{eq:unif_lim_inproof}. This yields~\eqref{eq:xi_small} and the proof is complete.
\end{proof}

\subsection{Step 3: Solutions of the characteristic equation}

Next we move towards proving existence of solutions of the characteristic equation~\eqref{eq:CharacteristicSign}. To this end, we will need the following auxiliary statement.

\begin{lemma}\label{lem:aux_eq}
Let $\nu=\ell+\frac{d}{2}-1$ and  $q=\frac{1-\varepsilon}{d-1}$ for a fixed $\varepsilon\in (0,1)$. Then there exists $h_0>0$ such that for all $h\geq h_0$, $\ell\in\mathcal{L}(h)$, and $j\in\mathcal{J}(h,\ell)$, the following claims hold:
\begin{enumerate}[(i)]
\item The function
\[
 \operatorname{err}_{\nu,j}(m):=-1+\frac{m}{4\pi j}\frac{\e^{\ii\theta_{\nu}(m)}}{\cos\theta_{\nu}(m)}\sqrt{1-\xi_{\nu}(m)}
\]
is analytic in $B_{2}(m_{\nu,j}^{(1)})$ and there is a~constant $C>0$ independent of $j$, $\ell$ and $m$ such that 
 \begin{equation}\label{eq:err_small}
  |\operatorname{err}_{\nu,j}(m)|\leq \frac{C}{\log^{q}h}
 \end{equation}
 for any $m\in B_{2}(m_{\nu,j}^{(1)})$.
\item If $m\in B_{2}(m_{\nu,j}^{(1)})$ satisfies
\begin{equation}
 \ii\left(\theta_{\nu}(m)-\frac{\pi}{4}-2\pi j\right)=\log\frac{4\pi j}{\sqrt{h}}+\log\left(1+\operatorname{err}_{\nu,j}(m)\right),
\label{eq:tow_char_eq}
\end{equation}
then $m$ is a solution of the characteristic equation~\eqref{eq:CharacteristicSign} with corresponding $k=k(m)$ given by \eqref{eq:defk}.
\end{enumerate}
\end{lemma}

\begin{proof}
\emph{Proof of claim~(i):} First, we verify the analyticity of the error function $\operatorname{err}_{\nu,j}$. Using Lemma~\ref{lem:xi_prop}, we may suppose $h_0$ to be sufficiently large to ensure that $\xi_{\nu}$ is analytic in $B_{2}(m_{\nu,j}^{(1)})$ and $|\xi_{\nu}(m)|<1$ for all $m\in B_{2}(m_{\nu,j}^{(1)})$, $\ell\in\mathcal{L}(h)$, $j\in\mathcal{J}(h,\ell)$, and $h\geq h_{0}$. Then the function $m\mapsto\sqrt{1-\xi_{\nu}(m)}$ is analytic in $B_{2}(m_{\nu,j}^{(1)})$. Since we already know from Lemma~\ref{lem:phase_analytic} and~\eqref{eq:B2_in_M} that $\theta_{\nu}$ is analytic in $B_{2}(m_{\nu,j}^{(1)})$, and also that $\cos\theta_{\nu}(m)\neq0$ in $B_{2}(m_{\nu,j}^{(1)})$, see the paragraph below~\eqref{eq:xi_inproof}, we conclude that $\operatorname{err}_{\nu,j}$ is indeed an analytic function in $B_{2}(m_{\nu,j}^{(1)})$.

Second, we prove~\eqref{eq:err_small}. We estimate the three factors in
\begin{equation}
  \operatorname{err}_{\nu,j}(m)=-1+\frac{m}{2\pi j}\frac{1}{1+\e^{-2\ii\theta_{\nu}(m)}}\sqrt{1-\xi_{\nu}(m)},
\label{eq:err_expre_inproof}
\end{equation}
using the uniform Landau symbol $\mathcal{O}$ for $h\to\infty$, where the involved constant is always independent of $\ell\in\mathcal{L}(h)$, $j\in\mathcal{J}(h,\ell)$, and $m\in B_{2}(m_{\nu,j}^{(1)})$. 
First, with the aid of~\eqref{eq:xi_small}, we deduce that 
\[
\sqrt{1-\xi_{\nu}(m)}=1+\bigO{\frac{1}{\log^{2q}h}}.
\]
Second, it follows from~\eqref{eq:im_theta_bounds} that
\[
 \frac{1}{1+\e^{-2\ii\theta_{\nu}(m)}}=\frac{1}{1+\bigO{h^{-2\alpha}}}=1+\bigO{h^{-2\alpha}}.
\]
Third, one infers from~\eqref{eq:re_m_reminder} and \eqref{eq:sup_im_lim} that
\[
 \frac{m}{2\pi j}=1+\bigO{\frac{1}{\log^{q}h}}.
\]
Inserting the last three estimates into~\eqref{eq:err_expre_inproof} amounts to the uniform asymptotic formula
\[
 \operatorname{err}_{\nu,j}(m)=\bigO{\frac{1}{\log^{q}h}},
\]
which yields~\eqref{eq:err_small}.

\emph{Proof of claim~(ii):} Using the definition~\eqref{eq:def_xi}, the  equation~\eqref{eq:Characteristic} can be written as
\[
 -\frac{\ii h}{m^{2}}=\frac{1-\xi_{\nu}(m)}{\cos^{2}\theta_{\nu}(m)}.
\]
Hence, any $m$ solving the equation
\[
 \e^{-\ii\left(\frac{\pi}{4}+2\pi j\right)}\frac{\sqrt{h}}{m}=\frac{\sqrt{1-\xi_{\nu}(m)}}{\cos\theta_{\nu}(m)}
\]
with an integer $j$, has to be a solution of~\eqref{eq:Characteristic}. When we rewrite the last equation as
\[
 \e^{\ii\left(\theta_{\nu}(m)-\frac{\pi}{4}-2\pi j\right)}=\frac{2m}{\sqrt{h}}\frac{\sqrt{1-\xi_{\nu}(m)}}{1+\e^{-2\ii\theta_{\nu}(m)}}=\frac{4\pi j}{\sqrt{h}}\left(1+\operatorname{err}_{\nu,j}(m)\right),
\]
we see that, if $m\in B_{2}(m_{\nu,j}^{(1)})$ is a solution of the equation~\eqref{eq:tow_char_eq}, then $m$ solves~\eqref{eq:Characteristic}. 
It is left to show that $m$ then also satisfies~\eqref{eq:CharacteristicSign}, by taking the principal branch square root on both sides of~\eqref{eq:Characteristic}. It suffices  to compare the signs of the leading order term on each side of~\eqref{eq:CharacteristicSign}.
By Lemma~\ref{lem:param_k}, we know that $k/m \to 1$. On the other hand, $\ii\tan\theta_\nu(m)\to 1$ by \eqref{eq:tan-cos_theta_bounds}, which, when inserted to~\eqref{eq:BesselProductAsymp}, implies that 
\[\frac{J_{\nu}'(m)H_{\nu}^{(1)}(k)}{J_{\nu}(m)(H_{\nu}^{(1)})'(k)}\to 1.\]
So the signs agree and the proof is complete.
\end{proof}

Claim~(ii) of Lemma~\ref{lem:aux_eq} lacks existence of the solutions. This is proven in the next statement.

\begin{Proposition}\label{prop:sol_m_j}
Let $\nu=\ell+\frac{d}{2}-1$. Then there exists $h_0>0$ such that for all $h\geq h_0$, $\ell\in\mathcal{L}(h)$, and $j\in\mathcal{J}(h,\ell)$, 
there exists a unique solution $m_{\nu,j}$ of the characteristic equation~\eqref{eq:CharacteristicSign} in the ball $B_{2}(m_{\nu,j}^{(1)})$. 
\end{Proposition}

\begin{proof}
The claim relies on Rouch{\' e}'s theorem and Lemma~\ref{lem:aux_eq}. Recalling definition~\eqref{eq:def_f_nu,j}, equation~\eqref{eq:tow_char_eq} can be written as 
\begin{equation}
 f_{\nu,j}(m)+\ii\log(1+\operatorname{err}_{\nu,j}(m))=0. 
\label{eq:tow_char_eq_inproof}
\end{equation}
We show that $h_{0}$ can be chosen so that, for all $h\geq h_{0}$, we have 
\begin{equation}
|f_{\nu,j}(m)|>|\log(1+\operatorname{err}_{\nu,j}(m))|
\label{eq:rouche_ineq_inproof}
\end{equation}
for all $m\in\partial B_{2}(m_{\nu,j}^{(1)})$, $\ell\in\mathcal{L}(h)$, and $j\in\mathcal{J}(h,\ell)$. Then, by Rouch{\' e}'s theorem, equation~\eqref{eq:tow_char_eq_inproof} has a unique solution $m_{\nu,j}$ in the ball $B_{2}(m_{\nu,j}^{(1)})$ since $m_{\nu,j}^{(1)}$ is a unique zero of $f_{\nu,j}$ in $B_{2}(m_{\nu,j}^{(1)})$, see Lemma~\ref{lem:aux_zeros}. According to claim~(ii) of Lemma~\ref{lem:aux_eq}, $m_{\nu,j}$ is then also a~solution of the characteristic equation~\eqref{eq:CharacteristicSign}.

Using definitions~\eqref{eq:def_psi}, 
\eqref{eq:def_f_nu,j} together with the equality $f_{\nu,j}(m_{\nu,j}^{(1)})=0$, we get
\[
 f_{\nu,j}(m)=f_{\nu,j}(m)-f_{\nu,j}(m_{\nu,j}^{(1)})=m-m_{\nu,j}^{(1)}-\psi_{\nu}(m)+\psi_{\nu}(m_{\nu,j}^{(1)}).
\]
Hence, if $m\in\partial B_{2}(m_{\nu,j}^{(1)})$, we find that
\[
 |f_{\nu,j}(m)|\geq 2-|\psi_{\nu}(m)-\psi_{\nu}(m_{\nu,j}^{(1)})|.
\]
Moreover, with the aid of the mean value theorem, the uniform limit~\eqref{eq:unif_lim_inproof}, and the estimate~\eqref{eq:psi_der_bound}, one may ensure that 
\[
|\psi_{\nu}(m)-\psi_{\nu}(m_{\nu,j}^{(1)})|<1
\]
for all $h$ large enough, $\ell\in\mathcal{L}(h)$, and $j\in\mathcal{J}(h,\ell)$. Thus,
$|f_{\nu,j}(m)|>1$ for all such $h$. On the other hand, by~\eqref{eq:err_small}, the right-hand side of~\eqref{eq:rouche_ineq_inproof} is smaller than $1$ for all $h$ sufficiently large. The inequality~\eqref{eq:rouche_ineq_inproof} follows.
\end{proof}

\begin{remark}
Notice that, given any $h\geq h_{0}$, $\ell\in\mathcal{L}(h)$, and $j_1,j_2\in\mathcal{J}(h,\ell)$, with $j_1\neq j_2$, we have $m_{\nu,j_1}\neq m_{\nu,j_2}$. This is a consequence of~\eqref{eq:B2_disjoint}.
\end{remark}

\subsection{Step 4: Inequalities for eigenvalues of $H_{h}$}

When estimating the left-hand side of~\eqref{eq:main-ineq}, we make use of the following inequalities for certain discrete eigenvalues of $H_{h}$.

\begin{Proposition}\label{prop:evls_ineq}
  Suppose $d\geq2$ and let  $\nu=\ell+\frac{d}{2}-1$. There exists $h_0>0$ such that, for all $h\geq h_0$, $\ell\in\mathcal{L}(h)$, and $j\in\mathcal{J}(h,\ell)$, the following claims hold true:
 \begin{enumerate}[(i)]
 \item The number
  \[
   \lambda_{\ell,j}:=\ii h+m_{\nu,j}^{2},
  \]
 where $m_{\nu,j}$ is the solution from Proposition~\ref{prop:sol_m_j},
 is an eigenvalue of $H_{h}$ of algebraic multiplicity at least~\eqref{eq:eigenspace}.
 \item We have the estimates
 \begin{equation}
  \Im\lambda_{\ell,j}\geq\frac{h}{2} 
  \quad\mbox{ and }\quad
  |\lambda_{\ell,j}|\leq (4\pi j)^{2}.
 \label{eq:lam_ineq}  
 \end{equation}
 \end{enumerate}
\end{Proposition}

\begin{proof}
\emph{Proof of claim~(i):}
 From the introductory analysis of the discrete spectrum of $H_{h}$ made in Section~\ref{sec:char_eq}, the claim~(i) follows if $m_{\nu,j}$ is a solution of the characteristic equation~\eqref{eq:CharacteristicSign} and the corresponding parameter $k_{\nu,j}$ given by the principal square root of $\ii h+m_{\nu,j}^{2}$ is of positive imaginary part, see~\eqref{eq:im_k_pos}. The former is established in Proposition~\ref{prop:sol_m_j}. The latter is verified next.
 
As the principal branch of the square root is used in the definition of $k_{\nu,j}$, we have $\Im k_{\nu,j}>0$ if 
\[
\Im\left(\ii h+m_{\nu,j}^{2}\right)=h+2\Re m_{\nu,j}\Im m_{\nu,j}>0.
\]
Since $m_{\nu,j}\in B_{2}(m_{\nu,j}^{(1)})$, see Proposition~\ref{prop:sol_m_j}, we may apply Lemma~\ref{lem:sup_re-im} and equation~\eqref{eq:set_J} to show that
\[
 \Re m_{\nu,j}=\bigO{j}=\bigO{h^{\gamma+1/2}}
 \quad\mbox{ and }\quad
 \Im m_{\nu,j}=\bigO{\log h},
\]
as $h\to\infty$, where the constants in the Landau $\mathcal{O}$s are independent of $\ell\in\mathcal{L}(h)$ and $j\in\mathcal{J}(h,\ell)$. Consequently, as $h\to\infty$, we have
\[
\Im\left(\ii h+m_{\nu,j}^{2}\right)=h\left[1+\bigO{h^{\gamma-1/2}\log h}\right]\!,
\]
from which we see that, by choosing $h_0$ large enough and using that $\gamma<1/2$, we have $\Im(\ii h+m_{\nu,j}^{2})>0$ for all $h\geq h_0$.

\emph{Proof of claim~(ii):} By the above computation, 
\[
\Im\lambda_{\ell,j}=h\left[1+\bigO{h^{\gamma-1/2}\log h}\right]\!,
\]
as $h\to\infty$, uniformly in $\ell\in\mathcal{L}(h)$ and $j\in\mathcal{J}(h,\ell)$. Therefore, for all $h$ sufficiently large, the first inequality in~\eqref{eq:lam_ineq} holds.

By using Lemma~\ref{lem:sup_re-im} again and the restrictions on the indices $\ell$ and $j$ from~\eqref{eq:restr_ell} and~\eqref{eq:set_J}, one checks that on the right-hand side of the inequality
\[
 \frac{|\lambda_{\ell,j}|}{(4\pi j)^{2}}\leq \frac{h+|m_{\nu,j}|^{2}}{(4\pi j)^{2}}=\frac{h}{(4\pi j)^{2}}+\left(\frac{\Re m_{\nu,j}}{4\pi j}\right)^{\!2}+\left(\frac{\Im m_{\nu,j}}{4\pi j}
\right)^{\!2}
\]
the first and third terms are of limit $0$, while the second term of limit $1/4$, as $h\to\infty$. Hence, we may suppose that $h_0$ is sufficiently large to guarantee that
\[
\frac{|\lambda_{\ell,j}|}{(4\pi j)^{2}}\leq 1
\]
for all $h\geq h_0$, which is the second inequality from~\eqref{eq:lam_ineq}. The proof is complete.
\end{proof}


Yet another auxiliary inequality will be needed.

\begin{lemma}\label{lem:multiplicityIneq}
Let $d\geq2$. For all sufficiently large $\ell\in\mathbb{N}$, we have
\begin{equation}\label{eq:multiplicityIneq}
	\binom{d+\ell-1}{d-1}-\binom{d+\ell-3}{d-1}\geq \frac{\ell^{d-2}}{(d-2)!}.
\end{equation}
\end{lemma}

\begin{proof}
 It is elementary to check, by using the definition of the binomial numbers, that the left-hand side of~\eqref{eq:multiplicityIneq} is a polynomial in $\ell$ of degree $d-2$ with the coefficient of $\ell^{d-2}$ equal to $2/(d-2)!$. This implies the statement.
\end{proof}

\subsection{Step 5: Proof of Proposition~\ref{LemmaafterMain}}

By Proposition~\ref{prop:evls_ineq} and Lemma~\ref{lem:multiplicityIneq}, there exists $h_0$ such that the inequalities from~\eqref{eq:lam_ineq} and~\eqref{eq:multiplicityIneq} hold for all $h\geq h_0$, $\ell\in\mathcal{L}(h)$, and $j\in\mathcal{J}(h,\ell)$. Therefore we may estimate 
\[
\sum_{\lambda\in\sigma_{\dd}(H_h)}\dfrac{(\im\lambda)^p}{|\lambda|^{d/2}}\geq \sum_{\ell\in\mathcal{L}(h)}\frac{\ell^{d-2}}{(d-2)!}\sum_{j\in\mathcal{J}(h,\ell)}\dfrac{(\im\lambda_{\ell,j})^p}{|\lambda_{\ell,j}|^{d/2}}
\geq C_{p,d}'\,h^{p}\sum_{\ell\in\mathcal{L}(h)}\ell^{d-2} \sum_{j\in\mathcal{J}(h,\ell)}\frac{1}{j^{d}}
\]
for all $h\geq h_{0}$, where 
\[
C_{p,d}':=\frac{1}{2^{p}(4\pi)^{d}(d-2)!}
\]
for now, but the notation $C_{p,d}'$ is used below for a generic constant depending only on $p$ and~$d$.
Thus, we have the lower bound
\begin{equation}
\frac{1}{h^{p}}\sum_{\lambda\in\sigma_{\dd}(H_h)}\dfrac{(\im\lambda)^p}{|\lambda|^{d/2}}\geq 
 C_{p,d}'\sum_{\ell\in\mathcal{L}(h)}\ell^{d-2} \sum_{j\in\mathcal{J}(h,\ell)}\frac{1}{j^{d}}
\label{eq:tow_hdk}
\end{equation}
for all $h$ sufficiently large.

Next we estimate from below the inner sum from the right-hand side of~\eqref{eq:tow_hdk}. Recalling definition~\eqref{eq:set_J}, we deduce
\[
\sum_{j\in\mathcal{J}(h,\ell)}\frac{1}{j^{d}}\geq\int_{1+\ell\log^{q}\ell}^{h^{\gamma+1/2}-1}\frac{\dd j}{j^{d}}=\frac{1}{d-1}\left[\frac{1}{(1+\ell\log^{q}\ell)^{d-1}}-\frac{1}{(h^{\gamma+1/2}-1)^{d-1}}\right].
\]
Since $\gamma>\beta>0$ the second term decays faster than the first one and therefore we may suppose, without loss of generality, that $h$ is large enough so that 
\[
\sum_{j\in\mathcal{J}(h,\ell)}\frac{1}{j^{d}}\geq 
\frac{1}{2(d-1)}\frac{1}{(\ell\log^{q}\ell)^{d-1}}=
\frac{1}{2(d-1)}\frac{1}{\ell^{d-1}\log^{1-\varepsilon}\ell},
\]
where definition~\eqref{eq:def_q} was used. Plugging the last estimate into~\eqref{eq:tow_hdk} yields
\begin{equation}
\frac{1}{h^{p}}\sum_{\lambda\in\sigma_{\dd}(H_h)}\dfrac{(\im\lambda)^p}{|\lambda|^{d/2}}\geq 
C_{p,d}'\sum_{\ell\in\mathcal{L}(h)}\frac{1}{\ell\log^{1-\varepsilon}\ell}
\label{eq:tow_hdk1}
\end{equation}
for all $h$ sufficiently large.

Similarly, for all $h$ sufficiently large, we may estimate the remaining sum by an integral as follows,
\[
\sum_{\ell\in\mathcal{L}(h)}\frac{1}{\ell\log^{1-\varepsilon}\ell}\geq\frac{1}{2}\int_{h^{\alpha+1/2}}^{h^{\beta+1/2}}\frac{\dd\ell}{\ell\log^{1-\varepsilon}\ell}
=\left[\left(\beta+\frac{1}{2}\right)^\varepsilon-\left(\alpha+\frac{1}{2}\right)^\varepsilon\right]\frac{\log^{\varepsilon}h}{2\varepsilon}.
\]
Taking also into account that the function 
\[
\varepsilon\mapsto\frac{1}{\varepsilon}\left[\left(\beta+\frac{1}{2}\right)^\varepsilon-\left(\alpha+\frac{1}{2}\right)^\varepsilon\right]
\]
is decreasing on $(0,1)$, we further deduce the lower bound with an $\varepsilon$-independent constant,
\[
\sum_{\ell\in\mathcal{L}(h)}\frac{1}{\ell\log^{1-\varepsilon}\ell}\geq\frac{\beta-\alpha}{2}\log^{\varepsilon}h.
\]
When the last estimate is used in~\eqref{eq:tow_hdk1}, we arrive at the inequality 
\[
\frac{1}{h^{p}}\sum_{\lambda\in\sigma_{\dd}(H_h)}\dfrac{(\im\lambda)^p}{|\lambda|^{d/2}}\geq 
 C_{p,d}'\,(\beta-\alpha)\log^{\varepsilon}h 
\]
for all $h$ sufficiently large. The proof of Proposition~\ref{LemmaafterMain}, and so of Theorem~\ref{ThmMainthm}, is complete.
\qed


\section*{Statements and Declarations}
\noindent \textbf{Data availability statement:} Data sharing not applicable to this article as no
datasets were generated or analysed during the current study.\\
\textbf{Conflict of interest declaration:} All authors certify that they have no affiliations
with or involvement in any organization or entity with any financial interest or nonfinancial
interest in the subject matter or materials discussed in this manuscript.\\
\textbf{Funding:} The research of F.~{\v S}. was supported by the GA{\v C}R grant No.\ 20-17749X. The PhD of S.~P.\ is funded through a Development and Promotion of Science and Technology (DPST) scholarship of the Royal Thai Government, ref.\ 5304.1/3758.

\bibliographystyle{acm}
 \bibliography{reference}

\end{document}